\documentstyle[epsf,12pt]{article}
\def\Bbb R{{\rm \bf R}}
\def\proclaim#1{\vskip2mm{\bf #1}\em}
\def\endproclaim{\em \vskip2mm}
\def\tag#1{\eqno(#1)}
\def\gathered{\begin{array}{c}}
\def\endgathered{\end{array}}
\def\text{\mbox}

\makeatletter
    
    \@addtoreset{figure}{section}
\makeatother

\begin{document}

\title {Two analytical formulae of the temperature inside a body by using partial lateral and initial data}
\author{Masaru IKEHATA\\
Department of Mathematics,
Graduate School of Engineering\\
Gunma University, Kiryu 376-8515, JAPAN}
\maketitle
\begin{abstract}
This paper considerers the problem of computing the value of a
solution of the heat equation at a given point inside a bounded
domain after the initial time.  It is assumed that the initial
value of the solution inside the domain (possibly in a part of the
domain) is known; the boundary value and the normal derivative on
a part of the boundary of the domain over a finite time interval
are known.  Two analytical formulae for the problem are given.
Both formulae make use of a special fundamental solution having a
large parameter of the backward heat equation.

\noindent
AMS: 35R30

\end{abstract}


\section{Introduction}
In this paper we consider the following inverse problem.  Assume
that we have a known heat conductive body and know the initial
temperature inside the body (say constant). We are in the
situation that we can not access the whole boundary, however, want
to know the time evolution of the temperature at a given point
inside the body (possibly close to an inaccessible part of the
boundary) after the initial time. How can one know it from the
temperature and heat flux on an accessible part of the boundary of
the body for a finite observation time?

This is a typical inverse and ill-posed problem and appears in many areas of 
engineering and medicine.
The aim of this paper is to develop an analytical approach to the problem.
We consider the simplest, however, important case:
the heat conductive body has a known isotropic and homogeneous heat
conductivity; there is no heat source or sink inside the body.
In this case, using Fourier's law, after a scaling, one may assume that the temperature
at a given point and time inside the body satisfies the heat
equation.  Then the problem is formulated as follows.

Let $\Omega\subset\Bbb R^n(n=1,2,3)$ be a bounded domain with a smooth boundary
and $0<T<\infty$.
Let $u=u(x,t)$ satisfy
$$\displaystyle
\partial_tu=\triangle u\,\,\,\text{in}\,\Omega\times]0,\,T[.
\tag {1.1}
$$
Given non empty open subsets $\Gamma\subset\partial\Omega\times]0,\,T[$ and $U\subset\Omega$
(typically $U=\Omega$)
find a formula for calculating $u(x,t)$ in $\Omega\times]0,\,T[$
from the data $u(x, 0)$ for $x\in U$ and $(u(x,t),\,\partial u/\partial\nu(x,t))$
for $(x,t)\in\Gamma\times\,]0,\,T[$.  Here $\nu$ denotes the unit
outward normal vector field to $\partial\Omega$.

We do not assume that $u(x,t)$ for
$(x,t)\in\,(\partial\Omega\times]0,\,T[)\setminus\Gamma$ and
$u(x,0)$ for $x\in\Omega\setminus U$ are {\it known}.  These are
considered as unknown input or output that one cannot control.
In this paper we give two solutions to the problem.
The consequence of the results stated in Sections 4 and 5 is the following.

Let $(x_0,t_0)\in\Omega\times]0,\,T[$, $\omega\in S^{n-1}$ and $c>0$.
If $(\Omega\setminus U)\times\{0\}$, $\Omega\times\{T\}$
and $(\partial\Omega\times\,]0,\,T[)\setminus\Gamma$ are contained in
the half space $(x\,\,t)^T\cdot\omega(c)\le (x_0\,\,t_0)^T\cdot\omega(c)-\delta$ in $\Bbb R^{n+1}$ with a $\delta>0$
and $\omega(c)=(1/\sqrt{1+c^2})(c\,\omega\,\,-1)^T(\in S^{n})$,
then there is an explicit formula to calculate $u(x_0,t_0)$ from the data.  See Figure 1.1 for the configuration
of $U$, $\Gamma$ and the half space $(x\,\,t)^T\cdot\omega(c)\le (x_0\,\,t_0)^T\cdot\omega(c)-\delta$.

\vspace{-1.5cm}

\begin{figure}[htbp]
\begin{center}
\epsfxsize=8.5cm
\epsfysize=12cm
\epsfbox{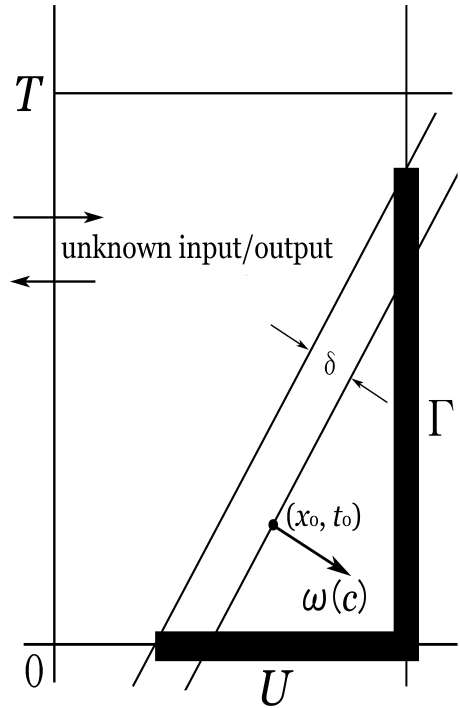}
\caption{Configuration.}\label{fig1}
\end{center}
\end{figure}

Here we explain our approach together with the construction of the paper.

\subsection{First approach}
First we observe that the complex exponential function $e^{x\cdot
z-t(z\cdot z)}$, $z\in\mbox{\boldmath $C$}^n$ satisfies the
backward heat equation $(\partial_t+\triangle)w=0$ in $\Bbb
R^{n+1}$.  Choosing a suitable $z$ depending on a large parameter
$\tau$, direction $\omega$ and $c$ (see (4.6)), we see that the
real part of the phase function $x\cdot z-t(z\cdot z)$ has the
form $\tau\sqrt{1+c^2}\,(x\,\,t)^T\cdot\omega(c)$.  This yields
that if $(x\,\,t)^T\cdot\omega(c)>0$, then $\vert e^{x\cdot
z-t(z\cdot z)}\vert\longrightarrow\infty$ as
$\tau\longrightarrow\infty$; if $(x\,\,t)^T\cdot\omega(c)<0$, then
$\vert e^{x\cdot z-t(z\cdot z)}\vert\longrightarrow 0$ as
$\tau\longrightarrow\infty$. This means that the asymptotic
behaviour of $e^{x\cdot z-t(z\cdot z)}$ divides the whole space
time into two parts whose common boundary is the hyper plane
$(x\,\,t)^T\cdot\omega(c)=0$.

Second we choose $D\subset\Bbb R^{n+1}$ with $\overline D\subset\Omega\times]0,\,T[$
and $(x_0,t_0)\in\partial D$ in such a way that for all $\rho$ smooth on $\overline D$ and for constants $\mu$
and $C$ independent of $\rho$
$$\displaystyle
\lim_{\tau\longrightarrow\infty}\tau^{\mu}e^{-x_0\cdot z+t_0(z\cdot z)}\int_D e^{x\cdot z-t(z\cdot z)}\rho(x,t)dxdt
=C\rho(x_0,t_0).
\tag {1.2}
$$
The examples of the choice of $D$ are given in subsections in Section 4.

Third we construct a special solution of the backward heat
equation with a special inhomogeneous term:
$$\displaystyle
\partial_tv+\triangle v+e^{x\cdot z-t(z\cdot z)}\chi_D(x,t)=0\,\,\text{in}\,\Bbb R^{n+1}.
\tag {1.3}
$$
The point of the property of $v$ is: the growth order of a suitable norm of
$e^{-x\cdot z+t(z\cdot z)}v$ as $\tau\longrightarrow\infty$ is at
most {\it algebraic}.
For the construction of $v$ we basically follow
a Fourier transform method for a reduced equation to (1.3) which has been used in \cite{SU}
since it is {\it constructive}.  It starts with introducing a special fundamental solution
denote by $G_z$ for the reduced equation.
This together with necessary estimates are described in Section 3.

Finally, integrating the equation (1.1) multiplied by $v$ over
$\Omega\times]0,\,T[$ and applying an integration by parts formula described in Section 2,
we know that the quantity
$$\displaystyle
\tau^{\mu}e^{-x_0\cdot z+t_0(z\cdot z)}\int_D e^{x\cdot z-t(z\cdot z)}u(x,t)dxdt
$$
which is coming from the inhomogeneous term in (1.3) multiplied by
$\tau^{\mu}e^{-x_0\cdot z+t_0(z\cdot z)}$,
is divided into two parts: the first part consists of the data on
$\Gamma$ and $U\times\{0\}$; the second unknown data on
$(\Omega\setminus U)\times\{0\}$,
$(\partial\Omega\times\,]0,\,T[)\setminus\Gamma$ and
$\Omega\times\{T\}$.  However, from the property of $v$ mentioned
above and the assumption on the configuration of $(\Omega\setminus
U)\times\{0\}$, $(\partial\Omega\times\,]0,\,T[)\setminus\Gamma$
and $\Omega\times\{T\}$ relative to the hyper plane
$(x\,\,t)^T\cdot\omega(c)=(x_0\,\,t_0)^T\cdot\omega(c)$ we know
that the second part tends to $0$ as $\tau\longrightarrow\infty$.
Combining this with (1.2) for $\rho=u$ on $\overline D$,
one gets a formula to calculate $u(x_0,t_0)$ from the data on
$\Gamma$ and $U\times\{0\}$ only.  See Section 4 for the precise description of the formula.
This approach can be considered as an
application of the {\it enclosure method} \cite{Ie3} to the heat equation (1.1).

\subsection{Second approach}
In the final section we give an integral representation formula of $G_z$ that is
important for calculating $v$.  Moreover as a byproduct we present a Carleman
type formula for the heat equation, which instead of
$v$, makes use of $e^{x\cdot z-t(z\cdot z)}G_z(x,t)$ directly
which is a special fundamental solution of the backward heat equation.
This approach is classical and the formula can be considered
as an extension to the heat equation of Yarmukhamedov's formula \cite{Y}
which has been established for the Laplace equation.

\subsection{Sideways heat equation}
There are extensive mathematical studies in the case when the one
dimensional heat equation is appropriate. This case is also
referred as the {\it sideways heat equation}. Carasso \cite{C}
considered the heat equation in the half line
$$\displaystyle
\partial_tu=u_{xx}\,\,\text{in}\,]0,\,\infty[\times\,]0,\,\infty[
$$
with the initial condition $u(x,0)=0$ for $0\le x<\infty$. He gave
a Tikhonov regularization procedure that yields an approximation
to $u(x,t)$ for $0<t<\infty$ at an arbitrary fixed $0<x<1$ from
noisy $u(1, t)$ for $0\le t<\infty$.  Note that by solving an
initial boundary value problem for the heat equation in
$1<x<\infty$, one gets the heat flux $u_x(1,t)$ from $u(1,t)$
and thus his problem is reduced to one-dimensional version
of our problem.
Motivated by Carasso's work Levine \cite{L} considered a radially
symmetric solution of the heat equation in higher dimensions and
established a Tikhonov regularization procedure. Their methods are
based on an explicit integral representation for the solution of
the heat equation and not time local in the sense: to determine
the value of the solution of the heat equation at a fixed point in
the body one needs the data for all the time ($T=\infty$). In
contrast to their methods our method does not make use of any
representation of the solution of the heat equation and needs data
only on an appropriate finite time interval depending on the
location of the point and time where and when we want to know the
value of the solution and covers fully multi dimensional cases.
For possible applications and numerical studies of the sideways
heat equation, see \cite{B, E} and references therein.

\section{A weak formulation of the direct problem and integration by parts formula}

In this section we describe what we mean by a solution of (1.1) together with an integration
by parts formula.  We follow the formulation of the direct problem described in \cite{DL}.

Given $\rho\in L^{\infty}(\partial\Omega)$, $f_0\in L^2(0,\,T;(H^1(\Omega))')$ and
$h_0\in L^2(0,\,T;H^{-1/2}(\partial\Omega))$ we say that $u\in W(0,\,T;H^1(\Omega), (H^1(\Omega))')$ satisfy
$$\begin{array}{c}
\displaystyle
\partial_tu-\triangle u=f_0\,\,\text{in}\,\Omega\times]0,\,T[,\\
\\
\displaystyle
\nabla u\cdot\nu+\rho u=h_0\,\,\text{on}\,\partial\Omega\times]0,\,T[
\end{array}
\tag {2.1}
$$
in the weak sense if the $u$ satisfies
$$\begin{array}{c}
\displaystyle
<u'(t),v>
+\int_{\Omega}\nabla u(x,t)\cdot\nabla v(x)dx
+\int_{\partial\Omega}u(t)\vert_{\partial\Omega}\cdot v\vert_{\partial\Omega}\,\rho dS\\
\\
\displaystyle
=<f_0(t),v>+<h_0(t),\,v\vert_{\partial\Omega}>\,\,\text{in}\,(0,\,T),
\end{array}
\tag {2.2}
$$
in the sense of distribution on $(0,\,T)$ for all $v\in H^1(\Omega)$.

Note that by Theorem 1 on p.473 in \cite{DL}
we see that every $u\in W(0,\,T;H^1(\Omega), (H^1(\Omega))')$ is almost everywhere equal to
a continuous function of $[0,\,T]$ in $L^2(\Omega)$.  Further, we have:
$$\displaystyle
W(0,\,T;H^1(\Omega), (H^1(\Omega))')\hookrightarrow C^0([0,\,T];L^2(\Omega)),
$$
the space $C^0([0,\,T];L^2(\Omega))$ being equipped with the norm
of uniform convergence. Thus one can consider $u(0)$ and $u(T)$ as
elements of $L^2(\Omega)$. Then by Theorems 1 and 2 on p.512 and
513 in \cite{DL} we see that given $u_0\in L^2(\Omega)$ there
exists a unique $u$ such that $u$ satisfies (2.1) in the weak
sense and satisfies the initial condition $u(0)=u_0$.

Given $\rho\in L^{\infty}(\partial\Omega)$, $f_1\in
L^2(0,\,T;(H^1(\Omega))')$ and $h_1\in
L^2(0,\,T;H^{-1/2}(\partial\Omega))$ we say that a $v\in
W(0,\,T;H^1(\Omega),(H^1(\Omega))')$ satisfy
$$\begin{array}{c}
\displaystyle
\displaystyle
\partial_tv+\triangle v=f_1\,\,\text{in}\,\Omega\times]0,\,T[,\\
\\
\displaystyle
\nabla v\cdot\nu+\rho v=h_1\,\,\text{on}\,\partial\Omega\times]0,\,T[
\end{array}
\tag {2.3}
$$
in the weak sense if the $v$ satisfies
$$\begin{array}{c}
\displaystyle
<v'(t),v>
-\int_{\Omega}\nabla u(x,t)\cdot\nabla v(x)dx
-\int_{\partial\Omega}u(t)\vert_{\partial\Omega}\cdot v\vert_{\partial\Omega}\,\rho dS\\
\\
\displaystyle
=<f_1(t),v>-<h_1(t),\,v\vert_{\partial\Omega}>\,\,\text{in}\,(0,\,T),
\end{array}
$$
in the sense of distribution on $(0,\,T)$ for all $v\in H^1(\Omega)$.

\proclaim{\noindent Proposition 2.1.} Let $u\in
W(0,\,T;H^1(\Omega), (H^1(\Omega))')$ satisfy (2.1) in a weak
sense; let $v\in W(0,\,T;H^1(\Omega), (H^1(\Omega))')$ satisfy
(2.3) in a weak sense.  Then the formula
$$\begin{array}{c}
\displaystyle
\int_0^T\left(<h_1(t),u(t)\vert_{\partial\Omega}>
-<h_0(t),v(t)\vert_{\partial\Omega}>\right)dt\\
\\
\displaystyle
=
\int_{\Omega}u(x,0)v(x,0)dx
-\int_{\Omega}u(x,T)v(x,T)dx
+\int_0^T(<f_1(t),u(t)>+<f_0(t),v(t)>)dt,
\end{array}
\tag {2.4}
$$
is valid.

\endproclaim

{\it\noindent Proof.}
Substituting $v=v(t)$ into (2.2), we have
$$\begin{array}{c}
\displaystyle
<u'(t),v(t)>+\int_{\Omega}\nabla u(x,t)\cdot\nabla v(x,t)dx\\
\\
\displaystyle
+\int_{\partial\Omega}u(t)\vert_{\partial\Omega}\cdot v(t)\vert_{\partial\Omega}\rho dS\\
\\
\displaystyle
=<f_0(t),v(t)>+<h_0(t),v(t)\vert_{\partial\Omega}>
\,\,\,{a.e.}\,\text{in}\,(0,\,T).
\end{array}
\tag {2.5}
$$
Similarly we have also
$$\begin{array}{c}
\displaystyle
<v'(t),u(t)>-\int_{\partial\Omega}\nabla v(x,t)\cdot\nabla u(x,t)dx\\
\\
\displaystyle
-\int_{\partial\Omega}v(t)\vert_{\partial\Omega}\cdot u(t)\vert_{\partial\Omega}\rho dS\\
\\
\displaystyle
=<f_1(t),u(t)>-<h_1(t),u(t)\vert_{\partial\Omega}>
\,\,\,{a.e.}\,\text{in}\,(0,\,T).
\end{array}
\tag {2.6}
$$
Taking the sum of (2.5) and (2.6), we have
$$\begin{array}{c}
\displaystyle
<u'(t),v(t)>+<v'(t),u(t)>
=<f_0(t),v(t)>+<h_0(t),v(t)\vert_{\partial\Omega}>\\
\\
\displaystyle
+<f_1(t),u(t)>-<h_1(t),u(t)\vert_{\partial\Omega}>.
\end{array}
$$
Integrating both sides of this equation over the interval $(0,\,T)$ and using the formula (Theorem 2 on p.477 in \cite{DL})
$$\begin{array}{c}
\displaystyle
\int_0^T<u'(t),v(t)>dt
+\int_0^T<v'(t),u(t)>dt\\
\\
\displaystyle
=\int_{\Omega}u(x,T)v(x,T)dx-\int_{\Omega}u(x,0)v(x,0)dx,
\end{array}
$$
we obtain (2.4).

\noindent
$\Box$

In particular, let $u\in W(0,\,T;H^1(\Omega), (H^1(\Omega))')$ satisfy
$$\begin{array}{c}
\displaystyle
\partial_tu-\triangle u=0\,\,\text{in}\,\Omega\times]0,\,T[,\\
\\
\displaystyle
\nabla u\cdot\nu+\rho u=h_0\,\,\text{on}\,\partial\Omega\times]0,\,T[
\end{array}
\tag {2.7}
$$
in the weak sense and $v\in W(0,\,T;H^1(\Omega),(H^1(\Omega))')$
satisfy
$$\begin{array}{c}
\displaystyle
\displaystyle
\partial_tv+\triangle v=f_1\,\,\text{in}\,\Omega\times]0,\,T[,\\
\\
\displaystyle
\nabla v\cdot\nu+\rho v=h_1\,\,\text{on}\,\partial\Omega\times]0,\,T[
\end{array}
\tag {2.8}
$$
in the weak sense.  Then (2.4) gives the formula
$$\begin{array}{c}
\displaystyle
\int_0^T\left(<h_1(t),u(t)\vert_{\partial\Omega}>
-<h_0(t),v(t)\vert_{\partial\Omega}>\right)dt\\
\\
\displaystyle
=\int_0^T<f_1(t),u(t)>dt
+\int_{\Omega}u(x,0)v(x,0)dx
-\int_{\Omega}u(x,T)v(x,T)dx.
\end{array}
\tag {2.9}
$$

\section{A solution of the equation
$\partial_t v+\triangle v+f=0$}

In this section given $f=f(x,t)$ we construct a special solution
of the equation
$$\displaystyle
-\partial_t v-\triangle v=f\,\,\text{in}\,\, \Bbb R^{n+1}.
\tag{3.1}
$$
Note that $-\partial_t-\triangle$ is the formal adjoint for
$\partial_t-\triangle$.

Given a complex vector $z\in\mbox{\boldmath $C$}^n$ set
$$
\displaystyle
v(x,t)=e^{x\cdot z-t(z\cdot z)}w(x,t),\,\,
f(x,t)=e^{x\cdot z-t(z\cdot z)}g(x,t).
$$
Then (3.1) becomes
$$\displaystyle
-\partial_t w-2z\cdot\nabla w-\triangle w=g.
\tag {3.2}
$$
Taking the Fourier transform of both sides, we obtain
$$\displaystyle
(-i\eta-2iz\cdot\xi+\vert\xi\vert^2)\hat{w}(\xi,\eta)
=\hat{g}(\xi,\eta).
$$

This motivates us to study the meaning of $1/P_z(\xi,\eta)$ where
$$\displaystyle
P_z(\xi,\eta)=-i\eta-2iz\cdot\xi+\vert\xi\vert^2.
$$
Let
$$\displaystyle
z=\mbox{\boldmath $a$}+i\mbox{\boldmath $b$},\,\, \mbox{\boldmath $a$}\,,\mbox{\boldmath $b$}\in\Bbb R^n.
$$
Since
$$
P_z(\xi,\eta)=\vert\xi+\mbox{\boldmath $b$}\vert^2-\vert\mbox{\boldmath $b$}\vert^2
-i(\eta+2\mbox{\boldmath $a$}\cdot\xi),
$$
we have:
$P_z(\xi,\eta)=0$ if and only if
$\displaystyle\vert\xi+\mbox{\boldmath $b$}\vert=\vert\mbox{\boldmath $b$}\vert$
and $\displaystyle\eta+2\mbox{\boldmath $a$}\cdot\xi=0$.
Therefore the set $\{(\xi,\eta)\in\Bbb R^{n+1}\,\vert\,P_z(\xi,\eta)=0\}$ is a compact set
and forms a submanifold of $\Bbb R^{n+1}$
with codimension $2$ provided $\mbox{\boldmath $b$}\not=0$.
This yields that for $z\in\mbox{\boldmath $C$}^n$ with $\text{Im}\,z\not=0$
the function $1/P_z(\xi,\eta)$
defines a tempered distribution on $\Bbb R^{n+1}$.

{\bf\noindent Definition 3.1.}
Define $G_z(x,t)$ as the inverse Fourier transform of $1/P_z(\xi,\eta)$:
$$
G_z(x,t)
=\frac{1}{(2\pi)^{n+1}}\int
e^{i(x\cdot\xi)+it\eta}\frac{d\xi d\eta}{ P_z(\xi,\eta)}
,\,\,\text{Im}\,z\not=0.
\tag {3.3}
$$

\noindent In this section we study the convolution operator
$g\mapsto G_z\ast g$ acting on the set of all rapidly decreasing
functions on $\Bbb R^{n+1}$.

First we restrict ourself to the case when
$\text{Re}\,z=0$.  More precisely, given $\omega\in S^{n-1}$ and $\tau>0$
set
$$\displaystyle
z=i\tau\omega.
$$
We set $F_{\tau}(x,t)= G_z(x,t)$ for this $z$, that is
$$
F_{\tau}(x,t)=
\frac{1}{(2\pi)^{n+1}}\int
e^{i(x\cdot\xi)+it\eta}\frac{d\xi d\eta}{ -i\eta+2\tau\omega\cdot\xi+\vert\xi\vert^2}.
$$

Let us study the property of $F_{\tau}(x,t)$.
For the purpose we employ the argument done in \cite{Ik0}.
The points are: a relationship between the operators $F_{\tau}\ast\,\cdot\,$ and $F_1\ast\,\cdot\,$;
an estimation of a scaling effect on weighted $L^2$-norms;
an weighted $L^2$-estimate for the operator $F_1\ast\,\cdot\,$.

\subsection{Scaling laws}

It is easy to see that $F_{\tau}(x,t)$ satisfies the scaling law:
$$\displaystyle
\forall\lambda>0\,\,
F_{\lambda\tau}(x,t)=\lambda^n F_{\tau}(\lambda x,\lambda^2 t).
\tag {3.4}
$$
Given a distribution $g(x,t)$ define
$$\displaystyle
g_{\lambda}(x,t)=g(\lambda x,\lambda^2 t ),\,\,\lambda>0.
$$
It is easy to see that $(g_{\lambda^{-1}})_{\lambda}=g$.
From (3.4) we have
$$\begin{array}{c}
\displaystyle
F_{\tau}\ast g(x,t)
=\int F_{\tau}(x-y,t-s)g(y,s) dyds\\
\\
\displaystyle
=\tau^n\int F_1(\tau(x-y), \tau^2(t-s))g(y,s)dyds.
\end{array}
$$
This yields
$$\begin{array}{c}
\displaystyle
(F_{\tau}\ast g)(\tau^{-1}x, \tau^{-2}t)
=\tau^n\int F_1(x-\tau y, t-\tau^2 s)g(y,s)dyds\\
\\
\displaystyle
=\frac{1}{\tau^{2}}
\int F_1(x-y, t-s)g(\tau^{-1}y, \tau^{-2}s)dyds.
\end{array}
$$
That is
$$\displaystyle
(F_{\tau}\ast g)_{\tau^{-1}}
=\frac{1}{\tau^{2}}F_1\ast(g_{\tau^{-1}})
$$
or equivalently
$$\displaystyle
F_{\tau}\ast g
=\frac{1}{\tau^{2}}\{F_1\ast (g_{\tau^{-1}})\}_{\tau}.
\tag {3.5}
$$
This yields also
$$\displaystyle
D_x^{\alpha}D_t^{\beta} F_{\tau}\ast g
=\tau^{-2+\vert\alpha\vert+2\beta}\left\{D_x^{\alpha}D_t^{\beta}F_1\ast(g_{\tau^{-1}})\right\}_{\tau}.
\tag {3.6}
$$

Given $s\in\Bbb R$ we denote by $L^2_s(\Bbb R^{n+1})$ the set of all tempered distributions $g=g(x,t)$ that
satisfies $(1+\vert x\vert^2+t^2)^{s/2}g\in L^2(\Bbb R^{n+1})$ and set
$$\displaystyle
\Vert g\Vert_s=\left(\int \vert g(x,t)\vert^2(1+\vert x\vert^2+t^2)^{s}dxdt\right)^{1/2},\,\,s\in\Bbb R.
$$
Note that the set of all rapidly decreasing functions on $\Bbb
R^{n+1}$ is dense in $L^2_s(\Bbb R^{n+1})$.

Given $R>0$ let $\tau\ge R$.
Set
$$\displaystyle
C(R)=\min\,\{R^4, R^2, 1\}\,(>0).
$$
Then we have
$$\begin{array}{c}
\displaystyle
1+\frac{\vert x\vert^2}{\tau^2}+\frac{t^2}{\tau^4}
=\frac{1}{\tau^4}(\tau^4+\tau^2\vert x\vert^2+t^2)\\
\\
\displaystyle
\ge\frac{1}{\tau^4}(R^4+R^2\vert x\vert^2+t^2)\\
\\
\displaystyle
\ge\frac{C(R)}{\tau^4}(1+\vert x\vert^2+t^2).
\end{array}
\tag {3.7}
$$
Since
$$\begin{array}{c}
\displaystyle
\Vert g_{\tau}\Vert^2_s
=\int \vert g(\tau x, \tau^2 t)\vert^2(1+\vert x\vert^2+t^2)^s dxdt\\
\\
\displaystyle
=\frac{1}{\tau^{n+2}}
\int\vert g(x,t)\vert^2\left(1+\frac{\vert x\vert^2}{\tau^2}+\frac{t^2}{\tau^4}\right)^s dxdt,
\end{array}
$$
from (3.7) one gets
$$\displaystyle
\Vert g_{\tau}\Vert_s
\le\frac{C(R)^{s/2}}{\tau^{2s+(n+2)/2}}\Vert g\Vert_s,\,\,s<0.
\tag {3.8}
$$

For $\tau\ge R$ we have
$$\begin{array}{c}
\displaystyle
1+\tau^2\vert x\vert^2+\tau^4t^2
=\tau^4(\tau^{-4}+\tau^{-2}\vert x\vert^2+t^2)\\
\\
\displaystyle
\le\tau^4\left(\frac{1}{R^4}+\frac{\vert x\vert^2}{R^2}+t^2\right)\\
\\
\displaystyle
\le \frac{\tau^4}{C(R)}(1+\vert x\vert^2+t^2).
\end{array}
\tag {3.9}
$$
Let $s'>0$.
Since
$$\displaystyle
\Vert g_{\tau^{-1}}\Vert^2_{s'}
=\tau^{n+2}\int \vert g(x,t)\vert^2(1+\tau^2\vert x\vert^2+\tau^4 t^2)^{s'}dxdt ,
$$
from (3.9) we obtain
$$\displaystyle
\Vert g_{\tau^{-1}}\Vert_{s'}
\le\frac{ \tau^{2s'+(n+2)/2}}{C(R)^{s'/2}}\Vert g\Vert_{s'}.
\tag {3.10}
$$

\subsection{Weighted $L^2$-estimates}

\proclaim{\noindent Lemma 3.1.} Let $-1<\delta<0$.  Given a rapidly
decreasing function $g$ on $\Bbb R^{n+1}$ the tempered
distribution $F_1\ast g$ belongs to $L^2_{\delta}(\Bbb R^{n+1})$
and there exists a positive constant $C_{\delta}$ independent of
$g$ and $\omega$ such that
$$\displaystyle
\Vert D_x^{\alpha}D_t^{\beta} F_1\ast g\Vert_{\delta}\le C_{\delta}\Vert g\Vert_{1+\delta},\,\vert\alpha\vert+2\beta\le 2.
\tag {3.11}
$$
\endproclaim

{\it\noindent Proof.}
For $z=i\omega$, we have $P_z(\xi,\eta)=\vert\xi+\omega\vert^2-1-i\eta$.
Let $\vert\xi\vert\ge 8$.
We have
$$\begin{array}{c}
\displaystyle
\vert P_z(\xi,\eta)\vert^2=(\vert\xi\vert^2+2\xi\cdot\omega)^2+\eta^2\\
\\
\displaystyle
=\vert\xi\vert^4+4(\xi\cdot\omega)\vert\xi\vert^2+4(\xi\cdot\omega)^2+\eta^2\\
\\
\displaystyle
\ge\vert\xi\vert^4+\eta^2-4\vert\xi\vert^3=\vert\xi\vert^4\left(1-\frac{4}{\vert\xi\vert}\right)+\eta^2\\
\\
\displaystyle
\ge\frac{1}{2}\vert\xi\vert^4+\eta^2
\ge\frac{1}{2}(\vert\xi\vert^4+\eta^2)
\ge\frac{1}{4}(\vert\xi\vert^2+\vert\eta\vert)^2.
\end{array}
$$
Next let $\vert\xi\vert\le 8$ and $\vert(\xi,\,\eta)\vert\ge 8\sqrt{1+8^2}$.
We have
$$
\displaystyle
\vert\eta\vert^2\ge 8^2(1+8^2)-\vert\xi\vert^2\ge 8^4\ge\vert\xi\vert^4.
$$
This yields $\vert\eta\vert\ge(\vert\xi\vert^2+\vert\eta\vert)/2$ and thus one gets
$$\displaystyle
\vert P_z(\xi,\eta)\vert^2\ge \eta^2\ge\left(\frac{1}{2}\right)^2(\vert\xi\vert^2+\vert\eta\vert)^2.
$$
Therefore it holds that, for all $(\xi,\eta)\in\Bbb R^{n+1}$ with $\vert(\xi,\,\eta)\vert\ge 8\sqrt{1+8^2}$
$$\displaystyle
\vert P_z(\xi,\eta)\vert\ge\frac{1}{2}(\vert\xi\vert^2+\vert\eta\vert).
$$
Using this inequality, a local representation of $1/P_z(\xi,\eta)$ in each neighbourhood of some zero points of $P_z(\xi,\eta)$
and Lemma 3.1 in \cite{SU}, we have the desired conclusion.

\noindent
$\Box$

A combination of (3.5), (3.6), (3.8), (3.10) and (3.11) yields

\proclaim{\noindent Proposition 3.1.} Let $-1<\delta<0$ and $R>0$.
For all rapidly decreasing functions $g$ on $\Bbb R^{n+1}$ and
$\tau\ge R$ we have
$$
\Vert D_{x}^{\alpha}D_t^{\beta}G_{i\tau\omega}\ast g\Vert_{\delta}
\le \frac{\displaystyle C_{\delta}\tau^{\vert\alpha\vert+2\beta}}{C(R)^{1/2}}\Vert g\Vert_{1+\delta},\,\vert\alpha\vert+2\beta\le 2.
\tag {3.12}
$$

\endproclaim

{\it\noindent Proof.}
First consider the case when $\vert\alpha\vert=\beta=0$.
We have
$$\begin{array}{c}
\displaystyle
\Vert F_{\tau}\ast g\Vert_{\delta}
=\frac{1}{\tau^{2}}\Vert\{F_1\ast(g_{\tau^{-1}})\}_{\tau}\Vert_{\delta}\\
\\
\displaystyle
\le \frac{C(R)^{\delta/2}}{\tau^{2\delta+2+(n+2)/2}}\Vert G_1\ast(g_{\tau^{-1}})\Vert_{\delta}\\
\\
\displaystyle
\le  \frac{C(R)^{\delta/2}C_{\delta}}{\tau^{2\delta+2+(n+2)/2}}
\Vert g_{\tau^{-1}}\Vert_{1+\delta}\\
\\
\displaystyle
\le \frac{C(R)^{\delta/2}C_{\delta}}{\tau^{2\delta+2+(n+2)/2}}
\frac{\tau^{2(1+\delta)+(n+2)/2}}{C(R)^{(1+\delta)/2}}\Vert g\Vert_{1+\delta}.
\end{array}
$$
Since we have (3.6), a similar argument yields (3.12) for the case when $\vert\alpha\vert\not=0$ or $\beta\not=0$.

\noindent
$\Box$

For general $z$ the following property is the starting point.

\proclaim{\noindent Proposition 3.2.}
For all $z\in\mbox{\boldmath $C$}^n$ with $\text{Im}\,z\not=0$ we have
$$\displaystyle
G_z(x,t)
=G_{\displaystyle i\text{Im}\,z}(x-2t\text{Re}\,z, t).
\tag {3.13}
$$

\endproclaim

{\it\noindent Proof.}
Let $z=\mbox{\boldmath $a$}+i\mbox{\boldmath $b$}$.  From (3.3) we have
$$\displaystyle
G_z(x,t)=\frac{1}{(2\pi)^{n+1}}\int
e^{i(x\cdot\xi)+it\eta}\frac{d\xi d\eta}
{-i(\eta+2\mbox{\boldmath $a$}\cdot\xi)+2\mbox{\boldmath $b$}\cdot\xi+\vert\xi\vert^2}.
$$
Then change of variables
$$\displaystyle
\eta'=\eta+2\mbox{\boldmath $a$}\cdot\xi,\,\,
\xi'=\xi
$$
yields
$$\displaystyle
G_z(x,t)=\frac{1}{(2\pi)^{n+1}}
\int e^{i (x\cdot\xi')+it(\eta'-2\mbox{\boldmath $a$}\cdot\xi')}
\frac{d\xi'd\eta'}
{-i\eta'-2i(i\mbox{\boldmath $b$})\cdot\xi'+\vert \xi'\vert^2}.
$$
Since
$$
i(x\cdot\xi')+it(\eta'-2\mbox{\boldmath $a$}\cdot\xi')
=i(x-2t\mbox{\boldmath $a$})\cdot\xi'+it\eta',
$$
we obtain the desired formula.

\noindent
$\Box$

Now given a real vector $\mbox{\boldmath $c$}$ set
$$\displaystyle
g_{\mbox{\boldmath $c$}}(x,t)=g(x-t\mbox{\boldmath $c$},t).
$$
From (3.13) we see that, for $z=\mbox{\boldmath $a$}+i\mbox{\boldmath $b$}$
$$\begin{array}{c}
\displaystyle
(G_{z}\ast g)(x,t)
=\int G_z(x-y,t-s)g(y,s)dyds\\
\\
\displaystyle
=\int G_{i\mbox{\boldmath $b$}}(x-y-2(t-s)\mbox{\boldmath $a$}, t-s)g(y,s)dyds.
\end{array}
$$
This yields that
$$\begin{array}{c}
\displaystyle
(G_z\ast g)(x+2t\mbox{\boldmath $a$},t)
=\int G_{i\mbox{\boldmath $b$}}(x-y+2s\mbox{\boldmath $a$}, t-s)g(y,s)dyds\\
\\
\displaystyle
=\int G_{i\mbox{\boldmath $b$}}(x-y',t-s)g(y'+2s\mbox{\boldmath $a$},s)dy'ds
\end{array}
$$
and thus we have
$$\displaystyle
(G_z\ast g)_{-2\mbox{\boldmath $a$}}
=G_{i\mbox{\boldmath $b$}}\ast(g_{-2\mbox{\boldmath $a$}}).
\tag {3.14}
$$
Since $(g_{\mbox{\boldmath $c$}})_{-\mbox{\boldmath $c$}}=g$, from (3.14) we obtain
$$\displaystyle
G_z\ast g=\{G_{i\mbox{\boldmath $b$}}\ast(g_{-2\mbox{\boldmath $a$}})\}_{2\mbox{\boldmath $a$}}
\tag {3.15}
$$
and also
$$\begin{array}{c}
\displaystyle
D_x^{\alpha}(G_z\ast g)=(D_x^{\alpha}G_{i\mbox{\boldmath $b$}}\ast g_{-2\mbox{\boldmath $a$}})_{2\mbox{\boldmath $a$}},\\
\\
\displaystyle
D_t(G_z\ast g)=-\sum_{j=1}^n
\left(\frac{\partial}{\partial x_j}G_{i\mbox{\boldmath $b$}}\ast g_{-2\mbox{\boldmath $a$}}\right)_{2\mbox{\boldmath $a$}}2a_j
+(D_tG_{i\mbox{\boldmath $b$}}\ast g_{-2\mbox{\boldmath $a$}})_{2\mbox{\boldmath $a$}}.
\end{array}
\tag {3.16}
$$

{\bf\noindent Remark 3.1.}
The equation (3.15) corresponds to the simple fact: a function $w(x,t)$ satisfies the equation (3.2)
if and only if the function $\tilde{w}(x,t)\equiv w(x+2t\mbox{\boldmath $a$},t)$ satisfies the equation
$$\displaystyle
-\partial_t\tilde{w}-2i\mbox{\boldmath $b$}\cdot\nabla\tilde{w}-\triangle\tilde{w}=g(x+2t\mbox{\boldmath $a$},t).
$$

Now we give an estimate for $G_z\ast g$.

\proclaim{\noindent Theorem 3.1.}
Let $-1<\delta<0$ and $R>0$.
Let $z=\mbox{\boldmath $a$}+i\mbox{\boldmath $b$}$.
For all rapidly decresing functions $g$ on $\Bbb R^{n+1}$ and
$\mbox{\boldmath $b$}\not=0$ with $\vert\mbox{\boldmath $b$}\vert\ge R$
we have
$$\begin{array}{c}
\displaystyle
\Vert D_x^{\alpha}G_{z}\ast g\Vert_{\delta}\le C(R,\delta)\left(\sqrt{1+\vert\mbox{\boldmath $a$}\vert^2}
+\vert\mbox{\boldmath $a$}\vert\right)\vert\mbox{\boldmath $b$}\vert^{\vert\alpha\vert}\Vert g\Vert_{\delta+1}, \,\vert\alpha\vert\le 2\\
\\
\displaystyle
\Vert D_t G_{z}\ast g\Vert_{\delta}\le C(R,\delta)\left(\sqrt{1+\vert\mbox{\boldmath $a$}\vert^2}
+\vert\mbox{\boldmath $a$}\vert\right)
(2\vert\mbox{\boldmath $a$}\vert\vert\mbox{\boldmath $b$}\vert+\vert\mbox{\boldmath $b$}\vert^2)\Vert g\Vert_{\delta+1}.
\end{array}
$$

\endproclaim

{\it\noindent Proof.}
It suffices to consider only the case when $\mbox{\boldmath $a$}\not=0$.
Let $\vert\alpha\vert=0$.
Set $f=G_{\displaystyle i\mbox{\boldmath $b$}}\ast(g_{\displaystyle -2\mbox{\boldmath $a$}})$.
From (3.15) we have
$$\begin{array}{c}
\displaystyle
\Vert G_z\ast g\Vert_{\delta}^2
=\Vert f_{\displaystyle 2\mbox{\boldmath $a$}}\Vert_{\delta}^2\\
\\
\displaystyle
=\int\vert f(x-2t\mbox{\boldmath $a$},t)\vert^2(1+\vert x\vert^2+t^2)^{\delta}dxdt\\
\\
\displaystyle
=\int\vert f(y,t)\vert^2(1+\vert y+2t\mbox{\boldmath $a$}\vert^2+t^2)^{\delta}dydt.
\end{array}
$$

One can write
$$\begin{array}{c}
\displaystyle
\vert y+2t\mbox{\boldmath $a$}\vert^2+t^2=A\left(\begin{array}{c} y\\
\\
\displaystyle
t\end{array}\right)
\cdot\left(\begin{array}{c} y\\
\\
\displaystyle
t
\end{array}\right)
\end{array}
$$
where
$$\displaystyle
A=\left(\begin{array}{cc}
\displaystyle I_n & 2\mbox{\boldmath $a$}\\
\\
\displaystyle
2\mbox{\boldmath $a$}^T & 1+4\vert\mbox{\boldmath $a$}\vert^2
\end{array}\right).
$$

It is easy to see that the eigenvalues $\lambda$ of $A$ coincides with the roots of the equation
$$\displaystyle
\lambda^2-2(1+2\vert\mbox{\boldmath $a$}\vert^2)\lambda+1=0.
$$
Solving this equation, we have
$$\displaystyle
\lambda=(1+2\vert\mbox{\boldmath $a$}\vert^2)\pm2\vert\mbox{\boldmath $a$}\vert
\sqrt{1+\vert\mbox{\boldmath $a$}\vert^2}.
$$

\noindent
Since the minimum eigenvalue has the form
$$\displaystyle
(1+2\vert\mbox{\boldmath $a$}\vert^2)-2\vert\mbox{\boldmath $a$}\vert
\sqrt{1+\vert\mbox{\boldmath $a$}\vert^2}
=\frac{1}
{\displaystyle (\sqrt{1+\vert\mbox{\boldmath $a$}\vert^2}+\vert\mbox{\boldmath $a$}\vert)^2},
$$
one has
$$\displaystyle
\frac{\vert y\vert^2+t^2}
{\displaystyle (\sqrt{1+\vert\mbox{\boldmath $a$}\vert^2}+\vert\mbox{\boldmath $a$}\vert)^2}
\le A\left(\begin{array}{c}
y\\
\\
t
\end{array}\right)
\cdot
\left(\begin{array}{c}
y\\
\\
t\end{array}
\right)
\le
(\sqrt{1+\vert\mbox{\boldmath $a$}\vert^2}+\vert\mbox{\boldmath $a$}\vert)^2
(\vert y\vert^2+t^2).
\tag {3.17}
$$

Since $-1<\delta<0$, from (3.17) we have
$$\displaystyle
\Vert f_{\displaystyle 2\mbox{\boldmath $a$}}\Vert_{\delta}^2
\le
\left(\sqrt{1+\vert\mbox{\boldmath $a$}\vert^2}+\vert\mbox{\boldmath $a$}\vert\right)^{-2\delta}
\Vert f\Vert_{\delta}^2
$$
and
$$
\Vert g_{\displaystyle -2\mbox{\boldmath $a$}}\Vert_{\delta+1}^2
\le
\left(\sqrt{1+\vert\mbox{\boldmath $a$}\vert^2}+\vert\mbox{\boldmath $a$}\vert\right)^{2(1+\delta)}
\Vert g\Vert_{\delta+1}^2.
$$
These together with (3.12) yield
$$\begin{array}{c}
\displaystyle
\Vert G_{z}\ast g\Vert_{\delta}^2
\le\left(\sqrt{1+\vert\mbox{\boldmath $a$}\vert^2}+\vert\mbox{\boldmath $a$}\vert\right)^{-2\delta}
\Vert f\Vert_{\delta}^2\\
\\
\displaystyle
\le C(R,\delta)\left(\sqrt{1+\vert\mbox{\boldmath $a$}\vert^2}+\vert\mbox{\boldmath $a$}\vert\right)^{-2\delta}
\Vert g_{\displaystyle -2\mbox{\boldmath $a$}}\Vert_{\delta+1}^2\\
\\
\displaystyle
\le C(R,\delta)\left(\sqrt{1+\vert\mbox{\boldmath $a$}\vert^2}+\vert\mbox{\boldmath $a$}\vert\right)^{-2\delta+2(1+\delta)}
\Vert g\Vert_{\delta+1}^2.
\end{array}
$$
Other cases also can be proved as above since we have (3.16).

\noindent
$\Box$

\noindent
Therefore the map $g\mapsto G_z\ast g\in L^2_{\delta}(\Bbb R^{n+1})$ can be uniquely extended
as a bounded linear operator of $L^2_{\delta+1}(\Bbb R^{n+1})$ into $L^2_{\delta}(\Bbb R^{n+1})$.
We denote it by the same symbol.
Then we see that, given $z$ with $\text{Im}\,z\not=0$ and
$g\in\,L^2_{\delta+1}(\Bbb R^{n+1})$ the the function
$$\displaystyle
v(x,t)=e^{x\cdot z-t(z\cdot z)}(G_{z}\ast g)(x,t)
\tag {3.18}
$$
satisfies the backward heat equation with a source term
in the sense of distribution:
$$\displaystyle
\partial_tv+\triangle v+e^{x\cdot z-t(z\cdot z)}g=0\,\,\text{in}\,\Bbb R^{n+1}.
\tag {3.19}
$$
Note that $e^{-x\cdot z+t(z\cdot z)}v(x,t)\in L^2_{\delta}(\Bbb R^{n+1})$
and this $v$ is {\it unique}.
This is a consequence of Theorem 7.1.27 of \cite{H}
and the facts that the set of all zero points of $P(\xi,\eta)$ has codimension 2 in $\Bbb R^{n+1}$
and $-1<\delta<0$.  See also Corollary 3.4 in \cite{SU} for this type of argument.

{\bf\noindent Remark 3.2.} Hsieh\cite{Hi} developed the scattering
theory associated with the operator $\partial_t-\triangle$ in
$\Bbb R^{2+1}$. For the purpose he studied the operator $L^1(\Bbb
R^{2+1})\ni f\mapsto (1/Q_z(\xi,\eta))\hat{q}\ast f\in L^1(\Bbb
R^{2+1})$ where $\hat{q}$ the Fourier transform  of a function $q$
on $\Bbb R^{2+1}$, $z\in\mbox{\boldmath $C$}^{2}$ and
$Q_z(\xi,\eta)$ the symbol of the operator $e^{-x\cdot z-t(z\cdot
z)}(\partial_t-\triangle)e^{x\cdot z+t(z\cdot z)}$. Note that
$e^{x\cdot z+t(z\cdot z)}$ satisfies the heat equation not the
backward heat equation. In the paper there is no result related
with Theorem 3.1.

\section{A computation formula of $u$ in $R^{n+1}$ with $n=1, 2, 3$}

In this section we assume that $h_0$ in (2.7) satisfies $h_0\in L^2(0,\,T;L^2(\partial\Omega))$
not just $h_0\in L^2(0,\,T;H^{-1/2}(\partial\Omega))$.

Let $u\in W(0,\,T;H^1(\Omega), (H^1(\Omega))')$
satisfy (2.7) in the weak sense.
Let $D$ a bounded open subset of $\Bbb R^{n+1}$ with $\overline D\subset\Omega\times\,]0,\,T[$.
We denote by $\chi_D$ the characteristic function of $D$.  Let $v$ be the distribution
given by (3.18) with $g=\chi_D$.

Now choose a sequence $g_j\in C^{\infty}_0(\Bbb R^{n+1})$ in such a way that
$g_j\longrightarrow \chi_D$ in $L^{2}_{\delta+1}(\Bbb R^{n+1})$ as $j\longrightarrow\infty$.
Define
$$
v_j(x,t)=e^{x\cdot z-t(z\cdot z)}(G_{z}\ast g_j)(x,t).
\tag {4.1}
$$
A combination of Theorem 3.1 and the Sobolev imbedding theorem in $\Bbb R^{n+1}$
we see that $v_j\in C^{\infty}(\Bbb R^{n+1})$ and $\{v_j\vert_{\Omega\times]0,\,T[}\}$
is Cauchy in $H^{2,1}(\Omega\times]0,\,T[)\equiv L^2(0,\,T;H^2(\Omega))
\cap H^1(0,\,T;L^2(\Omega))$(see pages 6-7 in \cite{LM}).
Since $v_j\vert_{\Omega\times]0,\,T[}\longrightarrow v\vert_{\Omega\times]0,\,T[}$
in $L^2(\Omega\times]0,\,T[)$, we conclude that $v\in H^{2,1}(\Omega\times]0,\,T[)$
and this $v$ satisfies
$$\displaystyle
\Vert D_x^{\alpha}D_t^{\beta}(e^{-x\cdot z+t(z\cdot z)}v(x,t))\Vert_{L^2(\Omega\times]0,\,T[)}
=O(\vert z\vert^3),\,\vert\alpha\vert+2\beta\le 2.
\tag {4.2}
$$
By the trace theorem (Theorem 2.1 on p. 9 in \cite{LM}) we have
$v_j\vert_{\partial\Omega\times]0,\,T[}\longrightarrow v\vert_{\partial\Omega\times]0,\,T[}$ in $H^{3/2,3/4}(\partial\Omega\times]0,\,T[)\equiv L^2(0,\,T;H^{3/2}(\partial\Omega))\cap H^{3/4}(0,\,T;L^2(\partial\Omega))$;
$\partial v_j/\partial\nu\vert_{\partial\Omega\times]0,\,T[}\longrightarrow\partial v/\partial\nu\vert_{\partial\Omega\times]0,\,T[}$
in $H^{1/2,1/4}(\partial\Omega\times]0,\,T[)\equiv L^2(0,\,T;H^{1/2}(\partial\Omega))\cap
H^{1/4}(0,\,T;L^2(\partial\Omega))$;
\newline{$v_j(x,0)\longrightarrow v(x,0)$ and $v_j(x,T)\longrightarrow v(x,T)$ in $H^1(\Omega)$.}

\noindent
Note that $v_j\in W(0,\,T;H^1(\Omega), (H^1(\Omega))')$
and satisfies (2.8) in the weak sense with $f_1=-e^{x\cdot z-t(z\cdot z)}g_j(x,t)$ and
$h_1=\partial v_j/\partial\nu+\rho v_j$
on $\partial\Omega\times]0,\,T[$.

Thus (2.9) yields
$$\begin{array}{c}
\displaystyle
\int_0^T\int_{\partial\Omega}\left\{\left(\frac{\partial v_j}{\partial\nu}+\rho v_j\right)u(t)\vert_{\partial\Omega}
-h_0(t)v_j(t)\vert_{\partial\Omega}\right\}dSdt\\
\\
\displaystyle
=-\int_0^T\int_{\Omega}e^{x\cdot z-t(z\cdot z)}g_j(x)u(x,t)dxdt
+\int_{\Omega}u(x,0)v_j(x,0)dx
-\int_{\Omega}u(x,T)v_j(x,T)dx.
\end{array}
$$
Taking the limit $j\longrightarrow\infty$, we obtain
$$\begin{array}{c}
\displaystyle
\int_{\partial\Omega\times]0,\,T[}\left\{\left(\frac{\partial v}{\partial\nu}(x,t)+\rho(x)v(x,t)\right)u(x,t)
-h_0(x,t)v(x,t)\right\}dSdt
\\
\\
\displaystyle
+\int_{\Omega}u(x,T)v(x,T)dx
-\int_{\Omega}u(x,0)v(x,0)dx\\
\\
\displaystyle
=-\int_D e^{x\cdot z-(z\cdot z)t}u(x,t)dxdt
\end{array}
\tag {4.3}
$$

\noindent Note that we made use of the fact that:  every $\phi\in
L^2(0,\,T;L^2(\partial\Omega))$ can be identified with
$\phi(x,t)\equiv \phi(t)(x)\in L^2(\partial\Omega\times]0,\,T[)$.

Divide $\partial\Omega\times]0,\,T[=\Gamma\cup\,\{(\partial\Omega\times]0,\,T[)\setminus\Gamma\}$
and $\Omega=U\cup(\Omega\setminus U)$.
Define
$$\begin{array}{c}
\displaystyle
I(\tau)
=
\int_{\Gamma}\left\{\left(\frac{\partial v}{\partial\nu}(x,t)+\rho(x)v(x,t)\right)u(x,t)-h_0(x,t)
v(x,t)\right\}dSdt\\
\\
\displaystyle
-\int_{U}v (x,0)u(x,0)dx.
\end{array}
$$
From (4.3) we have
$$\begin{array}{c}
\displaystyle
I(\tau)=-\int_D e^{x\cdot z-(z\cdot z)t}u(x,t)dxdt+R
\end{array}
\tag {4.4}
$$
where
$$\begin{array}{c}
\displaystyle
R=-\int_{\Omega}v(x,T)u(x,T)dx
+\int_{\Omega\setminus U}v(x,0)u(x,0)dx\\
\\
\displaystyle
-\int_{\partial\Omega\times]0,\,T[\setminus\Gamma}\left\{\left(\frac{\partial v}{\partial\nu}(x,t)+\rho(x)v(x,t)\right)u(x,t)-
h_0(x,t)v(x,t)\right\}dSdt.
\end{array}
\tag {4.5}
$$

Given $c>0$ and $\omega\in S^{n-1}$ define
$$\displaystyle
\omega(c)=\frac{1}{\sqrt{1+c^2}}\left(\begin{array}{c} c\,\omega\\
\\
\displaystyle
-1
\end{array}
\right)\in S^n
$$
and for $\tau$ with $c^2\tau>1$ set
$$z=
\left\{\begin{array}{lr}
\displaystyle
c\tau\left(\omega+i\sqrt{1-\frac{1}{c^2\tau}}\,\omega^{\perp}\right),\,\omega^{\perp}\in S^{n-1} & \quad\text{if $n=2,3$}\\
\\
\displaystyle
c\tau\left(1+i\sqrt{1-\frac{1}{c^2\tau}}\right)\omega,\,\omega=\{1,-1\} & \quad\text{if $n=1$.}
\end{array}
\right.
\tag {4.6}
$$
One can write
$$\displaystyle
\text{Re}\,\{x\cdot z-t(z\cdot z)\}
=\tau\sqrt{1+c^2}\left(\begin{array}{c} x\\
\\
t\end{array}
\right)
\cdot\omega(c)
$$
since $\text{Re}\,z\cdot z=\tau$.

In our method the concept introduced in the following plays an important role.

{\bf\noindent Definition 4.1.}
We say that $D$ is {\it visible} at $(x_0,t_0)\in\,\Bbb R^{n+1}$ as $\tau\longrightarrow\infty$
from the complex direction $z$ given by (4.6)
if there exist $\mu>0$ and constant $C\not=0$ such that
for all $\rho\in C^{\infty}(\overline D)$
$$\displaystyle
\lim_{\tau\longrightarrow\infty}e^{-x_0\cdot z+t_0(z\cdot z)}\tau^{\mu}\int_D e^{x\cdot z-t(z\cdot z)}\rho(x,t)dxdt
=C\rho(x_0,t_0).
\tag {4.7}
$$
The constant $C$ is unique if it exists.

\proclaim{\noindent Theorem 4.1.}
Let $(x_0,t_0)\in\Omega\times]0,\,T[$ be an arbitrary fixed point.
Assume that $T>0$, $\omega$, $\Gamma$ and $U$ satisfy the
following conditions:
$$\displaystyle
\sup_{x\in\Omega}\left(\begin{array}{c} x\\
\\
T\end{array}\right)\cdot\omega(c)
<\left(\begin{array}{c}
x_0\\
\\
t_0\end{array}\right)\cdot\omega(c);
\tag {4.8}
$$
$$\displaystyle
\sup_{x\in\Omega\setminus U}\left(\begin{array}{c} x\\
\\
0\end{array}\right)\cdot\omega(c)
<\left(\begin{array}{c}
x_0\\
\\
t_0\end{array}\right)\cdot\omega(c);
\tag {4.9}
$$
$$\displaystyle
\sup_{(x,t)\in(\partial\Omega\times]0,\,T[)\setminus\Gamma}\left(\begin{array}{c} x\\
\\
t\end{array}\right)\cdot\omega(c)
<\left(\begin{array}{c}
x_0\\
\\
t_0\end{array}\right)\cdot\omega(c).
\tag {4.10}
$$
Assume that
$D$ with $\overline D\subset\Omega\times]0,\,T[$ is visible at
$(x_0,t_0)$ from the complex direction $z$ given by (4.6).
Let $v$ be given by (3.18) with $g=\chi_D$.
Then we have
$$
\displaystyle
u(x_0,t_0)=-\frac{1}{C}
\lim_{\tau\longrightarrow\infty}
\tau^{\mu}e^{-x_0\cdot z+t_0(z\cdot z)} I(\tau),
\tag {4.11}
$$
where
$$\begin{array}{c}
\displaystyle
I(\tau)=
\int_{\Gamma}\left\{\left(\frac{\partial v}{\partial\nu}(x,t)+\rho(x)v(x,t)\right)u(x,t)
-h_0(x,t)v(x,t)\right\}dSdt
-\int_{U}v (x,0)u(x,0)dx.
\end{array}
$$

\endproclaim

{\it\noindent Proof.} Assumptions (4.8), (4.9), (4.10) together
with the trace theorem, (4.2) and (4.5) ensure that $\vert
\tau^{\mu}e^{-x_0\cdot z+t_0(z\cdot z)}R\vert$ is decaying as
$\tau\longrightarrow\infty$.  Since the heat operator
$\partial_t-\triangle$ is hypoelliptic, we know
$u\in\,C^{\infty}(\overline D)$. This together with (4.4) and
(4.7) yields (4.11).

\noindent
$\Box$

\noindent
The formula (4.11) can be considered as an application of an idea in \cite{Ie3, Ik3}
that was originally developed for the Cauchy problem for the stationary Schr\"odinger equation.

So the problem is reduced to: how to choose $D$ that is visible at $(x_0,t_0)$ from the complex direction $z$.
In the following subsections we consider this problem.

\subsection{The case $n=1$}

Let $\delta>0$.
We denote by $D(x_0,t_0,\omega(c),\delta)$ the inside of the triangle
with vertices $P=(x_0,t_0)$, $P_0=(x_0-(\delta/c)\sqrt{1+c^2}\omega,t_0)$
and $P_1=(x_0, t_0+\delta\sqrt{1+c^2}\,)$ in the space time $\Bbb R^{1+1}$.
The two points $P_1$ and $P_2$ are located on the line
$(x,t)^T\cdot\omega(c)=(x_0,t_0)^T\cdot\omega(c)-\delta$.

In \cite{Ik4} we have already known the following.  For the proof see
the proof of theorem 2.3 of in \cite{Ik4}.

\proclaim{\noindent Proposition 4.1.}
Let $D=D(x_0,t_0,\omega(c),\delta)$.  If $\rho\in C^{2}(\overline D)$,
then
$$\displaystyle
\lim_{\tau\longrightarrow\infty}
2(c\tau)^2\tau e^{-x_0\cdot z+t_0(z\cdot z)}\int_D e^{x\cdot z-t(z\cdot z)}\rho(x,t)dxdt
=
-\frac{i\vert P_1-P_0\vert^2\rho(P)}
{\displaystyle\vert P_1-P\vert\left(\sqrt{c^2+1}+i(c/\delta)\vert P_0-P\vert\right)}.
$$

\endproclaim

\noindent
Therefore this $D$ is visible at $(x_0,t_0)$ from
complex direction $z$.  The constants $\mu$ and $C$ in (4.7) are given by
$\mu=3$ and
$$\displaystyle
C
=-\frac{1+i}{4c^3}.
\tag {4.12}
$$

\subsection{The cases when $n=2,3$}

The cases when $n=2, 3$ start with describing the following which is easily derived
by the proof of Theorem 2.2 and Lemma 4.1 in \cite{Ik4}.

\proclaim{\noindent Proposition 4.2.} Let $n\ge 2$.  Let
$D\subset\Bbb R^{n+1}$ be a finite cone with a vertex at
$P=(x_0,t_0)$ and a bottom face $Q\not=\emptyset$ that is a
bounded open subset of $n$-dimensional hyper plane
$(x\,\,t)^T\cdot\omega(c)=(x_0\,\,t_0)^T\cdot\omega(c)-\delta$. If
$\rho\in C^{0,\theta}(\overline D)$ with $0<\theta\le 1$, then
$$\displaystyle
\lim_{\tau\longrightarrow\infty}\frac{2}{n!}
(c\tau)^{n+1}e^{-x_0\cdot z+t(z\cdot z)}
\int_D e^{x\cdot z-t(z\cdot z)}\rho(x,t)dxdt
=K_D\rho(P)
$$
where
$$\displaystyle
K_D=2\delta\int_Q
\frac{dS(y)}
{\displaystyle
\left(\frac{\delta\sqrt{c^2+1}}{c}-i(y-P)\cdot\left(\begin{array}{c}
\omega^{\perp}\\ 0
\end{array}
\right)\right)^{n+1}}.
$$
\endproclaim

Therefore if $K_D\not=0$, then constants $\mu$ and $C$ in (4.7)
are given by $\mu=n+1$ and
$$\displaystyle
C=\frac{n!K_D}{2c^{n+1}}.
\tag {4.13}
$$
However, it is not easy to show that
$K_D\not=0$ for $D$ with general $Q$.  In the following we specify $Q$ and
show that $K_D\not=0$.

\subsubsection{The case when $n=2$}

Let $\delta>0$.  Choose arbitrary two points $x_1$, $x_2$ on the line
$x\cdot\omega=x_0\cdot\omega-(\delta/c)\sqrt{1+c^2}$ in such a way that the orientation of the two vectors $\omega$, $x_1-x_2$
coincides with that of the standard basis $\mbox{\boldmath $e$}_1$, $\mbox{\boldmath $e$}_2$.
We denote by $D(x_0,x_1,x_2,\omega(c), \delta)$ the inside of the tetrahedron in $\Bbb R^{2+1}$
with the vertices $(x_0,t_0)$, $(x_1,t_0)$, $(x_2,t_0)$ and $(x_0,t_0+\delta\sqrt{1+c^2})$.
We see that the three points $(x_1, t_0)$, $(x_2,t_0)$ and $(x_0,t_0+\delta\sqrt{1+c^2})$
are on the plane $(x\,\,t)^T\cdot\omega(c)=(x_0,t_0)\cdot\omega(c)-\delta$.
Therefore $D(x_0,x_1,x_2,\omega(c),\delta)$ coincides with the finite cone with a vertex at $(x_0,t_0)$
and a bottom face $Q$ that is the triangle in $\Bbb R^{2+1}$ with the vertices $(x_1,t_0)$, $(x_2,t_0)$ and
$(x_0,t_0+\delta\sqrt{1+c^2})$.

From (4.2) in \cite{Ik4} we obtain
the formula
$$\begin{array}{c}
\displaystyle
K_D
=K_D\vartheta\cdot(0\,\,-1)^T\\
\\
\displaystyle
=
c^3\sum_{j=1}^3
\frac{\vert(\mbox{\boldmath $\nu$}_j\times\mbox{\boldmath $\nu$}_{j-1})
\times(\mbox{\boldmath $\nu$}_{j+1}\times\mbox{\boldmath $\nu$}_{j})\vert}
{\{(\mbox{\boldmath $\nu$}_j\times\mbox{\boldmath $\nu$}_{j-1})
\cdot\vartheta)\}
\{(\mbox{\boldmath $\nu$}_{j+1}\times\mbox{\boldmath $\nu$}_{j})
\cdot\vartheta)\}}
\mbox{\boldmath $\nu$}_j\cdot(0\,\,-1)^T
\end{array}
\tag {4.14}
$$
where $D=D(x_0,x_1,x_2,\omega(c),\delta)$, $\mbox{\boldmath $\nu$}_1=\mbox{\boldmath $\nu$}_4$,
$\mbox{\boldmath $\nu$}_2$, $\mbox{\boldmath $\nu$}_3=\mbox{\boldmath $\nu$}_0$ are
the unit outward normal vector to the faces of $D$ that are triangles $\Delta_1$
with the vertices $(x_0, t_0), (x_1, t_0), (x_0,t_0+\delta\sqrt{1+c^2})$, $\Delta_2$
with the vertices $(x_0,t_0), (x_2, t_0), (x_0,t_0+\delta\sqrt{1+c^2})$, $\Delta_3$ with the vertices
$(x_0,t_0), (x_1,t_0), (x_2,t_0)$; $\vartheta=(c(\omega+i\omega^{\perp})\,\,-1)^T$.

Since Corollary 4.1 in \cite{Ik4} ensures this $K_D\not=0$, we conclude that $D$
is visible at $(x_0,t_0)$ from the complex direction $z$.  Note that
$\mbox{\boldmath $\nu$}_3=(0\,\,-1)^T$ and $\mbox{\boldmath $\nu$}_1\cdot\mbox{\boldmath $\nu$}_3=
\mbox{\boldmath $\nu$}_2\cdot\mbox{\boldmath $\nu$}_3=0$.  Therefore from (4.13) and (4.14) we have the simpler
expression
$$\displaystyle
C=\frac{\vert(\mbox{\boldmath $\nu$}_3\times\mbox{\boldmath $\nu$}_{2})
\times(\mbox{\boldmath $\nu$}_{1}\times\mbox{\boldmath $\nu$}_{3})\vert}
{\{(\mbox{\boldmath $\nu$}_3\times\mbox{\boldmath $\nu$}_{2})
\cdot\vartheta)\}
\{(\mbox{\boldmath $\nu$}_{1}\times\mbox{\boldmath $\nu$}_{3})
\cdot\vartheta)\}}.
\tag {4.15}
$$

Now choose $x_1$, $x_2$, $\delta$ in such a way that $\overline{D(x_0,x_1,x_2,\omega(c),\delta)}\subset\Omega\times]0,\,T[$.
Then we obtain the formula (4.11) for $D=D(x_0,x_1,x_2,\omega(c),\delta)$, $\mu=3$ and $C$ given by (4.15).

\subsubsection{The case when $n=3$}

Let $\delta>0$.  Choose arbitrary three points $x_1$, $x_2$ and
$x_3$ on the plane
$x\cdot\omega=x_0\cdot\omega-(\delta/c)\sqrt{1+c^2}$ in such a way
that the orientation of the three vectors $\omega$, $x_1-x_2$,
$x_3-x_2$ coincides with that of the standard basis
$\mbox{\boldmath $e$}_1$, $\mbox{\boldmath $e$}_2$,
$\mbox{\boldmath $e$}_3$.  The four points $x_0$, $x_1$, $x_2$ and $x_3$ form
a tetrahedron $\Delta$ in $\Bbb R^3$.  We denote by $\mbox{\boldmath $\nu$}$
the unit outward normal vector field to $\partial\Delta$.  $\partial\Delta$
consists of four triangles: $T_1$ with the vertices $x_0$, $x_1$ and $x_2$;
$T_2$ with the vertices $x_0$, $x_3$ and $x_2$; $T_3$ with the vertices $x_0$, $x_1$ and $x_3$;
$T_4$ with the vertices $x_1$, $x_2$ and $x_3$.  Since $\mbox{\boldmath $\nu$}$
takes a constant vector on each $T_j$, we denote the vector by $\mbox{\boldmath $\nu$}_j$.
In particular, we have $\mbox{\boldmath $\nu$}_4=-\omega$.

We denote by $D(x_0,x_1,x_2, x_3,
\omega(c), \delta)$ the inside of the finite cone in $\Bbb
R^{3+1}$ with a vertex $(x_0,t_0+\delta\sqrt{1+c^2})$ and the
bottom that is the inside of the tetrahedron in the space $t=t_0$
with vertices $(x_0,t_0)$, $(x_1, t_0)$, $(x_2, t_0)$ and $(x_3,
t_0)$.  Then the boundary of $D(x_0,x_1,x_2, x_3, \omega(c),
\delta)$ consists of five tetrahedrons: $\Delta_1$ with the
vertices $(x_0,t_0)$, $(x_1, t_0)$, $(x_2,t_0)$ and $(x_0,
t_0+\delta\sqrt{1+c^2})$; $\Delta_2$ with the vertices
$(x_0,t_0)$, $(x_2,t_0)$, $(x_3, t_0)$ and $(x_0,t_0+\delta\sqrt{1+c^2})$;
$\Delta_3$ with the vertices
$(x_0,t_0)$, $(x_3,t_0)$, $(x_1,t_0)$ and
$(x_0,t_0+\delta\sqrt{1+c^2})$;
$\Delta_4$ with the vertices $(x_0,t_0)$, $(x_1,t_0)$, $(x_2,t_0)$
and $(x_3,t_0)$; $Q$ with vertices $(x_1,t_0)$,
$(x_2,t_0)$, $(x_3,t_0)$ and $(x_0, t_0+\delta\sqrt{1+c^2})$.

We see that $D=D(x_0,x_1,x_2,x_3,\omega(c),\delta)$ coincides with
the finite cone with a vertex at $(x_0,t_0)$ and the bottom $Q$.
Let $\mbox{\boldmath $a$}$ be an arbitrary constant vector in
$\mbox{\boldmath $C$}^{3+1}$. Since
$$\displaystyle
\nabla_{(x,t)}\cdot(e^{x\cdot z-t(z\cdot z)}\mbox{\boldmath $a$})
=(z\,\,-\tau)^T\cdot\mbox{\boldmath $a$}\,e^{x\cdot z-t(z\cdot z)},
$$
we have
$$\begin{array}{c}
\displaystyle
(z\,\,-\tau)^T\cdot\mbox{\boldmath $a$}\int_D e^{x\cdot z-t(z\cdot z)}dxdt
=\int_D\nabla_{(x,t)}\cdot(e^{x\cdot z-t(z\cdot z)}\mbox{\boldmath $a$})dxdt\\
\\
\displaystyle
=\sum_{j=1}^3\mbox{\boldmath $a$}\cdot(\mbox{\boldmath $\nu$}_j\,\,0)^T
\int_{\Delta_j} e^{x\cdot z-t(z\cdot z)}dS(x,t)
+\mbox{\boldmath $a$}\cdot(0\,\,-1)^T
\int_{\Delta_4} e^{x\cdot z-t(z\cdot z)}dS(x,t)\\
\\
\displaystyle
-\mbox{\boldmath $a$}\cdot\omega(c)\int_Q e^{x\cdot z-t(z\cdot z)}dS(x,t),
\end{array}
$$
where we made use of the fact that the unit outward normal vector
to $\partial D$ takes $(\mbox{\boldmath $\nu$}_j\,\,0)^T$ on
$\Delta_j$ for each $j=1,2,3$; $-\mbox{\boldmath
$e$}_4=(0\,\,-1)^T$ on $\Delta_4$; $-\omega(c)$ on $Q$. Since
$\mbox{\boldmath $a$}$ is arbitrary, one obtains
$$\begin{array}{c}
\displaystyle
\left(\begin{array}{c}
z\\
-\tau
\end{array}
\right)\int_D e^{x\cdot z-t(z\cdot z)}dxdt\\
\\
\displaystyle
=\sum_{j=1}^3
\left(\begin{array}{c}
\mbox{\boldmath $\nu$}_j\\
0\end{array}
\right)
\int_{\Delta_j} e^{x\cdot z-t(z\cdot z)}dS(x,t)
+\left(\begin{array}{c}
0\\
-1\end{array}
\right)
\int_{\Delta_4} e^{x\cdot z-t(z\cdot z)}dS(x,t)\\
\\
\displaystyle
-\omega(c)\int_Q e^{x\cdot z-t(z\cdot z)}dS(x,t).
\end{array}
\tag {4.16}
$$

Since $Q$ is
included in the hyper plane
$(x\,\,t)^T\cdot\omega(c)=(x_0\,\,t_0)^T\cdot\omega(c)-\delta$,
we have
$$\displaystyle
e^{-x_0\cdot z+t_0(z\cdot z)}\omega(c)\int_Q e^{x\cdot z-t(z\cdot z)}dS(x,t)
=O(e^{-\tau\delta\sqrt{1+c^2}}).
\tag {4.17}
$$

We compute the integral
$$\displaystyle
I_j(\tau)=\int_{\Delta_j}e^{x\cdot z-t(z\cdot z)}dS(x,t),\,\,j=1,2,3,4.
$$

\noindent
On $\Delta_j$ for each $j$ one can write
$$\displaystyle
\left(\begin{array}{c}
x\\
t
\end{array}
\right)
=
\left(\begin{array}{c} x_0\\
t_0\end{array}\right)
+\alpha\mbox{\boldmath $a$}_j+\beta\mbox{\boldmath $b$}_j+\gamma\mbox{\boldmath $c$}_j,
$$
where $(\alpha,\beta,\gamma)\in\Delta_0
=\{(\alpha,\beta,\gamma)\,\vert\alpha+\beta+\gamma\le
1,\,\,\alpha,\beta,\gamma\ge 0\}$ and $\mbox{\boldmath $a$}_j$,
$\mbox{\boldmath $b$}_j$, $\mbox{\boldmath $c$}_j$ are suitable
linearly independent vectors in $\Bbb R^{3+1}$ and satisfy the
condition
$$\displaystyle
\mbox{\boldmath $a$}_j\cdot\omega(c)<0,\,\mbox{\boldmath $b$}_j\cdot\omega(c)<0,\,\,
\mbox{\boldmath $c$}_j\cdot\omega(c)<0.
\tag {4.18}
$$

Writing
$A_j=(\mbox{\boldmath $a$}_j\,\,\mbox{\boldmath $b$}_j\,\,\mbox{\boldmath $c$}_j)$ which is
a $4\times 3$-matrix and $\tau\Delta_0=\{(\alpha,\beta,\gamma)\,\vert \alpha+\beta+\gamma\le\tau,\,\,
\alpha,\beta,\gamma\ge 0\}$, we have
$$\begin{array}{c}
\displaystyle
e^{-x_0\cdot z+t_0(z\cdot z)}I_j(\tau)=\sqrt{\text{det}\,A_j^TA_j}
\int_{\Delta_0} e^{(\alpha\mbox{\boldmath $a$}_j+\beta\mbox{\boldmath $b$}_j+\gamma\mbox{\boldmath $c$}_j)
\cdot(z\,\,-\tau)^T} d\alpha d\beta d\gamma\\
\\
\displaystyle
=\frac{1}{\tau^3}
\int_{\tau\Delta_0}
e^{(\alpha\mbox{\boldmath $a$}_j+\beta\mbox{\boldmath $b$}_j+\gamma\mbox{\boldmath $c$}_j)
\cdot(c(\omega+i\sqrt{1-(1/c^2\tau)}\omega^{\perp})\,\,-1)^T}
d\alpha d\beta d\gamma.
\end{array}
$$
Thus together with (4.18) yields
$$\begin{array}{c}
\displaystyle
\lim_{\tau\longrightarrow\infty}\tau^3 e^{-x_0\cdot z+t_0(z\cdot z)}I_j(\tau)
=
\int_0^{\infty}d\alpha\int_0^{\infty}d\beta\int_0^{\infty}d\gamma
e^{(\alpha\mbox{\boldmath $a$}_j+\beta\mbox{\boldmath $b$}_j+\gamma\mbox{\boldmath $c$}_j)
\cdot\mbox{\boldmath $\vartheta$}}\\
\\
\displaystyle
=\frac{(-1)^3}{(\mbox{\boldmath $a$}_j\cdot\mbox{\boldmath $\vartheta$})
(\mbox{\boldmath $b$}_j\cdot\mbox{\boldmath $\vartheta$})
(\mbox{\boldmath $c$}_j\cdot\mbox{\boldmath $\vartheta$})},
\end{array}
\tag {4.19}
$$
where $\mbox{\boldmath $\vartheta$}=(c(\omega+i\omega^{\perp})\,\,-1)^T$.

From (4.16), (4.17), (4.19) we obtain
$$\begin{array}{c}
\displaystyle
\lim_{\tau\longrightarrow\infty}\left(\begin{array}{c}
z\\
-\tau\end{array}
\right)\tau^3 e^{-x_0\cdot z+t_0(z\cdot z)}
\int_De^{x\cdot z-t(z\cdot z)}dxdt\\
\\
\displaystyle
=-\sum_{j=1}^3
\frac{\sqrt{\mbox{det}\,A_j^TA_j}}
{(\mbox{\boldmath $a$}_j\cdot\mbox{\boldmath $\vartheta$})
(\mbox{\boldmath $b$}_j\cdot\mbox{\boldmath $\vartheta$})
(\mbox{\boldmath $c$}_j\cdot\mbox{\boldmath $\vartheta$})}
\,\left(\begin{array}{c}
\mbox{\boldmath $\nu$}_j\\
0\end{array}\right)
-\frac{\sqrt{\mbox{det}\,A_4^TA_4}}
{(\mbox{\boldmath $a$}_4\cdot\mbox{\boldmath $\vartheta$})
(\mbox{\boldmath $b$}_4\cdot\mbox{\boldmath $\vartheta$})
(\mbox{\boldmath $c$}_4\cdot\mbox{\boldmath $\vartheta$})}
\left(\begin{array}{c}
0\\
-1\end{array}
\right).
\end{array}
$$
This yields the formula
$$
\displaystyle
\lim_{\tau\longrightarrow\infty}
e^{-x_0\cdot z+t_0(z\cdot z)}\tau^4
\int_De^{x\cdot z-t(z\cdot z)}dxdt
=-\frac{\sqrt{\mbox{det}\,A_4^TA_4}}
{(\mbox{\boldmath $a$}_4\cdot\mbox{\boldmath $\vartheta$})
(\mbox{\boldmath $b$}_4\cdot\mbox{\boldmath $\vartheta$})
(\mbox{\boldmath $c$}_4\cdot\mbox{\boldmath $\vartheta$})}.
$$
By choosing $\rho\equiv 1$ in Proposition 4.2,
one concludes
$$\displaystyle
K_D=-\frac{2}{3!}c^4\frac{\sqrt{\mbox{det}\,A_4^TA_4}}
{(\mbox{\boldmath $a$}_4\cdot\mbox{\boldmath $\vartheta$})
(\mbox{\boldmath $b$}_4\cdot\mbox{\boldmath $\vartheta$})
(\mbox{\boldmath $c$}_4\cdot\mbox{\boldmath $\vartheta$})}.
$$
Therefore $D$ is visible at $(x_0,t_0)$ from the complex direction $z$
and (4.7) is valid for $\mu=4$ and $C$ given by the formula
$$\displaystyle
C=
-\frac{\sqrt{\mbox{det}\,A_4^TA_4}}
{(\mbox{\boldmath $a$}_4\cdot\mbox{\boldmath $\vartheta$})
(\mbox{\boldmath $b$}_4\cdot\mbox{\boldmath $\vartheta$})
(\mbox{\boldmath $c$}_4\cdot\mbox{\boldmath $\vartheta$})}.
\tag {4.20}
$$
Therefore (4.11) is valid for this $C$ and $\mu=4$ provided $\delta$ is chosen
in such a way that $\overline{D(x_0,x_1,x_2,x_3,\omega(c),\delta)}\subset \Omega\times\,]0,\,T[$.

\section{An integral representation of $G_z$ and a byproduct}

It is quite important for us to compute $v$ given by (3.18) with $g=\chi_D$.
In this section we give an integral representation of the distribution $G_z(x,t)$ together with
$$\displaystyle
K_z(x,t)=e^{x\cdot z-t(z\cdot z)}G_{z}(x,t),\,\,z=\mbox{\boldmath $a$}+i\mbox{\boldmath $b$},\, \mbox{\boldmath $b$}\not=0
\tag {5.1}
$$
which is a solution of the equation $\partial_t v+\triangle v+\delta(x,t)=0$ in $\Bbb R^{n+1}$.
Using the representation of $K_z$, we show that $K_z$ in the hyper space $(x\,\,t)^T\cdot\omega(c)<0$
with $z$ given by (4.6) is exponentially decaying as $\tau\longrightarrow\infty$.  As a byproduct of this fact
we see that $K_z$ yields a Carleman type formula for the heat equation.

\subsection{Representation of $G_z$}

\proclaim{\noindent Proposition 5.1.}
It holds that
$$\begin{array}{c}
\displaystyle
G_{z}(x,t)\\
\\
\displaystyle
=e^{\displaystyle -i(x-2t\mbox{\boldmath $a$})\cdot\mbox{\boldmath $b$}-\vert\mbox{\boldmath $b$}\vert^2 t}\times\\
\\
\displaystyle
\left\{-\left(\frac{\vert\mbox{\boldmath $b$}\vert}{2\pi}\right)^n
\int_{\vert\xi\vert<1}e^{\displaystyle i\vert\mbox{\boldmath $b$}\vert (x-2t\mbox{\boldmath $a$})\cdot\xi}
e^{\displaystyle\vert\xi\vert^2\vert\mbox{\boldmath $b$}\vert^2t}d\xi
+H(-t)\left(\frac{1}{2\sqrt{\pi\vert t\vert}}\right)^n
e^{\displaystyle\frac{\vert x-2t\mbox{\boldmath $a$}\vert^2}{4t}}\right\}.
\end{array}
\tag {5.2}
$$
\endproclaim

{\it\noindent Proof.}
First we give a representation for $F_1(x,t)=G_{i\omega}(x,t)$
where $\omega=\mbox{\boldmath $b$}/\vert\mbox{\boldmath $b$}\vert$.
The starting point is the following formulae
for the Fourier transform:

Let $H(t)$ denote the Heaviside function.
Let $\text{Re}\,c\not=0$.  It is easy to see that:

\noindent
if $\text{Re}\,c>0$, then
$$\displaystyle
\frac{1}{-i\eta-c}
=-\int H(t)e^{-ct}e^{-it\eta}dt;
$$

\noindent
if $\text{Re}\,c<0$, then
$$\displaystyle
\frac{1}{-i\eta-c}
=\int H(-t)e^{-ct}e^{-it\eta}dt.
$$

\noindent
Thus $F_1(x,t)$ becomes
$$\begin{array}{c}
\displaystyle
F_1(x,t)=
\frac{1}{(2\pi)^{(n+1)}}\int
e^{i(x\cdot\xi)+it\eta}\frac{d\xi d\eta}{ -i\eta-(1-\vert\xi+\omega\vert^2)}\\
\\
\displaystyle
=\frac{1}{(2\pi)^{n}}
\{\int_{\vert\xi+\omega\vert<1}e^{ix\cdot\xi}
\{-H(t)e^{(\vert\xi+\omega\vert^2-1)t}\}d\xi
+\int_{\vert\xi+\omega\vert>1}e^{ix\cdot\xi}
\{H(-t)e^{(\vert\xi+\omega\vert^2-1)t}\}d\xi\}\\
\\
\displaystyle
=-\frac{H(t)}{(2\pi)^n}\int_{\vert\xi+\omega\vert<1}e^{ix\cdot\xi}
e^{(\vert\xi+\omega\vert^2-1)t}d\xi
+\frac{H(-t)}{(2\pi)^n}\int_{\vert\xi+\omega\vert>1}e^{ix\cdot\xi}
e^{(\vert\xi+\omega\vert^2-1)t}d\xi.
\end{array}
$$
This together with a change of variables yields
$$\begin{array}{c}
\displaystyle
F_1(x,t)\\
\\
\displaystyle
=-\frac{H(t)e^{-ix\cdot\omega}}{(2\pi)^n}
\int_{\vert\xi\vert<1}e^{ix\cdot\xi}
e^{(\vert\xi\vert^2-1)t}d\xi
+\frac{H(-t)e^{-ix\cdot\omega}}{(2\pi)^n}\int_{\vert\xi\vert>1}e^{ix\cdot\xi}
e^{(\vert\xi\vert^2-1)t}d\xi
\end{array}
$$
and thus we obtain
$$\begin{array}{c}
\displaystyle
e^{ix\cdot\omega+t}F_1(x,t)\\
\\
\displaystyle
=-\frac{H(t)}{(2\pi)^n}
\int_{\vert\xi\vert<1}e^{ix\cdot\xi}
e^{\vert\xi\vert^2t}d\xi
+\frac{H(-t)}{(2\pi)^n}\int_{\vert\xi\vert>1}e^{ix\cdot\xi}
e^{\vert\xi\vert^2t}d\xi.
\end{array}
\tag {5.3}
$$

Since
$$\displaystyle
\int e^{ix\cdot\xi}e^{-\vert\xi\vert^2a}d\xi
=\left(\sqrt{\frac{\pi}{a}}\right)^ne^{-\frac{\vert x\vert^2}{4a}},\,\,a>0,
\tag {5.4}
$$
we have, for all $t<0$
$$\begin{array}{c}
\displaystyle
\frac{1}{(2\pi)^n}\int_{\vert\xi\vert>1}e^{ix\cdot\xi}e^{\vert\xi\vert^2t}d\xi
=\frac{1}{(2\pi)^n}
\left\{\left(\sqrt{\frac{\pi}{\vert t\vert}}\right)^ne^{\frac{\vert x\vert^2}{4t}}
-\int_{\vert\xi\vert<1}e^{ix\cdot\xi}e^{\vert\xi\vert^2t}d\xi\right\}\\
\\
\displaystyle
=\left(\frac{1}{2\sqrt{\pi\vert t\vert}}\right)^ne^{\frac{\vert x\vert^2}{4t}}
-\frac{1}{(2\pi)^n}\int_{\vert\xi\vert<1}e^{ix\cdot\xi}
e^{\vert\xi\vert^2t}d\xi.
\end{array}
$$
This together with (5.3) yields
$$\begin{array}{c}
\displaystyle
e^{ix\cdot\omega+t}F_1(x,t)\\
\\
\displaystyle
=-\frac{1}{(2\pi)^n}\{H(t)+H(-t)\}\int_{\vert\xi\vert<1}e^{ix\cdot\xi}e^{\vert\xi\vert^2t}d\xi
+H(-t)\left(\frac{1}{2\sqrt{\pi\vert t\vert}}\right)^ne^{\frac{\vert x \vert^2}{4t}}\\
\\
\displaystyle
=-\frac{1}{(2\pi)^n}\int_{\vert\xi\vert<1}e^{ix\cdot\xi}e^{\vert\xi\vert^2t}d\xi
+H(-t)\left(\frac{1}{2\sqrt{\pi\vert t\vert}}\right)^n e^{\displaystyle\frac{\vert x\vert^2}{4t}}.
\end{array}
\tag {5.5}
$$

It follows from (3.3) that
$$\displaystyle
G_{\displaystyle i\mbox{\boldmath $b$}}(x,t)=\vert\mbox{\boldmath $b$}\vert^n F_1(\vert\mbox{\boldmath $b$}\vert x,
\vert\mbox{\boldmath $b$}\vert^2 t).
$$

\noindent
A combination of this and (5.5) gives
$$\begin{array}{c}
\displaystyle
G_{\displaystyle i\mbox{\boldmath $b$}}(x,t)\\
\\
\displaystyle
=e^{\displaystyle -ix\cdot\mbox{\boldmath $b$}-\vert\mbox{\boldmath $b$}\vert^2 t}
\left\{-\left(\frac{\vert\mbox{\boldmath $b$}\vert}{2\pi}\right)^n
\int_{\vert\xi\vert<1}e^{\displaystyle i\vert\mbox{\boldmath $b$}\vert x\cdot\xi}
e^{\displaystyle\vert\xi\vert^2\vert\mbox{\boldmath $b$}\vert^2t}d\xi
+H(-t)\left(\frac{1}{2\sqrt{\pi\vert t\vert}}\right)^n
e^{\displaystyle\frac{\vert x\vert^2}{4t}}\right\}.
\end{array}
$$
From this we immediately obtain (5.2) since we have (3.13).

\noindent
$\Box$

\subsection{Representation of $K_z$}

From (5.1), (5.2) and the equation
$$\displaystyle
x\cdot\mbox{\boldmath $a$}-t\vert\mbox{\boldmath $a$}\vert^2+\frac{1}{4t}\vert x-2t\mbox{\boldmath $a$}\vert^2
=\frac{\vert x\vert^2}{4t},
\tag {5.6}
$$
it follows that
$$\displaystyle
K_z(x,t)=
H(-t)\left(\frac{1}{2\sqrt{\pi\vert t\vert}}\right)^ne^{\displaystyle\frac{\vert x\vert^2}{4t}}
+w_z(x,t)
\tag {5.7}
$$
where
$$\displaystyle
w_z(x,t)=-e^{\displaystyle x\cdot\mbox{\boldmath $a$}-t\vert\mbox{\boldmath $a$}\vert^2}\left(\frac{\vert\mbox{\boldmath $b$}\vert}{2\pi}\right)^n
\int_{\vert\xi\vert<1}e^{\displaystyle i\vert\mbox{\boldmath $b$}\vert(x-2t\mbox{\boldmath $a$})\cdot\xi}\,
e^{\displaystyle\vert\mbox{\boldmath $b$}\vert^2\vert\xi\vert^2 t}d\xi.
\tag {5.8}
$$

Remarks are in order.

(i)  Since the distribution
$$\displaystyle
H(-t)\left(\frac{1}{2\sqrt{\pi\vert t\vert}}\right)^ne^{\displaystyle\frac{\vert x\vert^2}{4t}}
$$
is the fundamental solution of the equation $\partial_t v+\triangle v=0$, the $w_z(x,t)$ is an entire solution
of the backward heat equation.  Moreover $w_z$
is a smooth function on the whole space.

(ii)  From (5.8) and a change of variables we know that $w_z(x,t)=\overline{w_z(x,t)}$.  Thus $w_z$ is real valued
and hence we have
$$\displaystyle
w_z(x,t)=-e^{\displaystyle x\cdot\mbox{\boldmath $a$}-t\vert\mbox{\boldmath $a$}\vert^2}\left(\frac{\vert\mbox{\boldmath $b$}\vert}{2\pi}\right)^n
\int_{\vert\xi\vert<1}\cos\,(\vert\mbox{\boldmath $b$}\vert(x-2t\mbox{\boldmath $a$})\cdot\xi)\,
e^{\displaystyle\vert\mbox{\boldmath $b$}\vert^2\vert\xi\vert^2 t}d\xi.
$$
Note also that $w_z$ does not depend on the direction of the vector $\mbox{\boldmath $b$}$.

(iii)  Write
$$\displaystyle
\vert\mbox{\boldmath $b$}\vert^2\vert\xi\vert^2t+i\vert\mbox{\boldmath $b$}\vert(x-2t\mbox{\boldmath $a$})\cdot\xi
=t\left(\vert\mbox{\boldmath $b$}\vert\xi+i\frac{1}{2t}(x-2t\mbox{\boldmath $a$})\right)^2
+\frac{1}{4t}\vert x-2t\mbox{\boldmath $a$}\vert^2.
$$
Combining this with (5.6), we can rewrite (5.8) as
$$\displaystyle
w_z(x,t)=
-\frac{\displaystyle e^{\vert x\vert^2/(4t)}}{(2\pi)^n}
\int_{\vert\xi\vert<\vert\mbox{\boldmath $b$}\vert}
e^{\displaystyle t(\xi+i(x-2t\mbox{\boldmath $a$})/(2t))^2}d\xi.
$$

We study more the expression (5.7).

(i)  The case when $t>0$.
Since
$$\displaystyle
x\cdot\mbox{\boldmath $a$}-t\vert\mbox{\boldmath $a$}\vert^2+t\vert\mbox{\boldmath $b$}\vert^2
=\tau\sqrt{1+c^2}\left(\begin{array}{c} x\\
t\end{array}\right)\cdot\omega(c),
\tag {5.9}
$$
it follows from (5.7) and (5.8) that
$$\begin{array}{c}
\displaystyle
K_z(x,t)
=-e^{\tau\sqrt{1+c^2}\,(x\,\,t)^T\cdot\omega(c)}
\left(\frac{\vert\mbox{\boldmath $b$}\vert}{2\pi}\right)^n
\int_{\vert\xi\vert<1}
e^{\displaystyle i\vert\mbox{\boldmath $b$}\vert(x-2t\mbox{\boldmath $a$})\cdot\xi}
e^{\displaystyle -(1-\vert\xi\vert^2)\vert\mbox{\boldmath $b$}\vert^2 t}d\xi.
\end{array}
\tag {5.10}
$$

(ii) The case when $t<0$.
It follows from (5.4) that
$$\begin{array}{c}
\displaystyle
\int_{\vert\xi\vert<1}
e^{\displaystyle i\vert\mbox{\boldmath $b$}\vert(x-2t\mbox{\boldmath $a$})\cdot\xi}
e^{\displaystyle\vert\mbox{\boldmath $b$}\vert^2\vert\xi\vert^2 t}d\xi\\
\\
\displaystyle
=\frac{1}{\vert\mbox{\boldmath $b$}\vert^n}
\left(\sqrt{\frac{\pi}{-t}}\,\right)^n
e^{\displaystyle\vert x-2t\mbox{\boldmath $a$}\vert^2/(4t)}
-\int_{\vert\xi\vert>1}
e^{\displaystyle i\vert\mbox{\boldmath $b$}\vert(x-2t\mbox{\boldmath $a$})\cdot\xi}
e^{\displaystyle\vert\mbox{\boldmath $b$}\vert^2\vert\xi\vert^2 t}d\xi.
\end{array}
$$
This together with (5.7) and (5.8) yields
$$\displaystyle
K_z(x,t)
=e^{\tau\sqrt{1+c^2}\,(x\,\,t)^T\cdot\omega(c)}
\left(\frac{\vert\mbox{\boldmath $b$}\vert}{2\pi}\right)^n
\int_{\vert\xi\vert>1}
e^{\displaystyle i\vert\mbox{\boldmath $b$}\vert(x-2t\mbox{\boldmath $a$})\cdot\xi}
e^{\displaystyle\vert\mbox{\boldmath $b$}\vert^2(\vert\xi\vert^2-1)t}d\xi.
\tag {5.11}
$$

{\bf\noindent Remark 5.1.}
Using the well known formula
$$\displaystyle
\int_{\vert\xi\vert=r}e^{i\eta\cdot\xi}dS(\xi)
=(2\pi)^{n/2}r^{n/2}\vert\eta\vert^{-(n-2)/2}J_{(n-2)/2}(\vert\eta\vert r),\,\,\forall\eta\in\Bbb R^n,
$$
one can rewrite (5.8), (5.10) and (5.11) as one-dimensional integrals.

\subsection{Exponential decaying of $K_z$ in the hyper space $(x\,\,t)^T\cdot\omega(c)<0$
and a Carleman type formula for the heat equation}

In this subsection first we show that $K_z(x,t)$ is exponentially decaying as $\tau\longrightarrow\infty$
if $(x\,\,t)^T\cdot\omega(c)<0$ and $z$ is given by (4.6).

\proclaim{\noindent Proposition 5.2.}
Given $\delta>0$ we have, as $\tau\longrightarrow\infty$
$$\displaystyle
\sup_{\displaystyle (x\,\,t)^T\cdot\omega(c)<-\delta}\vert K_z(x,t)\vert=O(e^{-\tau\sqrt{1+c^2}\delta} \tau^n).
$$

\endproclaim

{\it\noindent Proof.}
Let $(x,\,t)$ satisfy $\displaystyle (x\,\,t)^T\cdot\omega(c)<-\delta$.

\noindent
(i)  The case when $t>0$.  From (5.10) we have immediately
$$\displaystyle
\vert K_z(x,t)\vert\le C_ne^{-\tau\sqrt{1+c^2}\delta}\vert\mbox{\boldmath $b$}\vert^n.
\tag {5.12}
$$

\noindent (ii) The case when $t<0$.  We divide this case into two
subcases: (a)  $\vert\mbox{\boldmath $b$}\vert^2 (-t)$ is large;
(b) $\vert\mbox{\boldmath $b$}\vert^2 (-t)$ is not large and can
be arbitrary small.

First consider (a).  Given $R>0$ let $t$ satisfy $\vert\mbox{\boldmath $b$}\vert^2 t<-R$.
From (5.4) and (5.11) we obtain
$$\displaystyle
\vert K_z(x,t)\vert\le e^{-\tau\sqrt{1+c^2}\,\delta}\left(\frac{\vert\mbox{\boldmath $b$}\vert}{2\pi}\right)^n
\int_{\vert\xi\vert>1}e^{-R(\vert\xi\vert^2-1)}d\xi
\le C_ne^{-\tau\sqrt{1+c^2}\,\delta}\vert\mbox{\boldmath $b$}\vert^n\frac{e^R}{R^{n/2}}.
\tag {5.13}
$$

Next consider (b).  We employ the expression (5.7) and (5.8).
Let $t$ satisfy $-R\le \vert\mbox{\boldmath $b$}\vert^2 t<0$.
Using (5.9), we can rewrite (5.8) as
$$\displaystyle
w_z(x,t)=-e^{\displaystyle \tau\sqrt{1+c^2}\,(x\,\,t)^T\cdot\omega(c)}e^{\displaystyle
-\vert\mbox{\boldmath $b$}\vert^2t}
\left(\frac{\vert\mbox{\boldmath $b$}\vert}{2\pi}\right)^n
\int_{\vert\xi\vert<1}e^{\displaystyle i\vert\mbox{\boldmath $b$}\vert(x-2t\mbox{\boldmath $a$})\cdot\xi}\,
e^{\displaystyle\vert\mbox{\boldmath $b$}\vert^2\vert\xi\vert^2 t}d\xi.
\tag {5.14}
$$
Since $t<0$ and $-\vert\mbox{\boldmath $b$}\vert^2 t\le R$, it follows from (5.14) that
$$\displaystyle
\vert w_z(x,t)\vert\le C_n e^{-\tau\sqrt{1+c^2}\delta}\vert\mbox{\boldmath $b$}\vert^n e^{R}.
\tag {5.15}
$$

Since $cx\cdot\omega<cx\cdot\omega-t<-\delta\sqrt{1+c^2}$, we have
$\vert x\cdot\omega\vert>(\delta/c)\sqrt{1+c^2}$ and thus $\vert x\vert>(\delta/c)\sqrt{1+c^2}$.
Using this together with $\vert\mbox{\boldmath $b$}\vert^2\vert t\vert\le R$,
we obtain
$$\begin{array}{c}
\displaystyle
\vert t\vert^{-n}e^{\vert x\vert^2/(4t)}
=(\vert\mbox{\boldmath $b$}\vert^2\vert t\vert)^{-n}
\vert\mbox{\boldmath $b$}\vert^{2n}
e^{\displaystyle -(\vert\mbox{\boldmath $b$}\vert^2\vert x\vert^2)/(4\vert\mbox{\boldmath $b$}\vert^2\vert t\vert)}\\
\\
\displaystyle
\le\frac{\vert\mbox{\boldmath $b$}\vert^{2n}}
{R^n}e^{\displaystyle -(\vert\mbox{\boldmath $b$}\vert^2(\delta/c)^2(1+c^2))/(4R)}.
\end{array}
$$
This together with (5.7) and (5.15) yields that
$$\displaystyle
\vert K_z(x,t)\vert
\le C_n\left(e^{-\tau\sqrt{1+c^2}\delta}\vert\mbox{\boldmath $b$}\vert^n e^R
+\frac{\vert\mbox{\boldmath $b$}\vert^{2n}}
{R^n}e^{\displaystyle -(\vert\mbox{\boldmath $b$}\vert^2(\delta/c)^2(1+c^2))/(4R)}\right).
\tag {5.16}
$$
A combination of (5.13) and (5.16) gives
$$\displaystyle
\vert K_z(x,t)\vert
\le C_n\left\{e^{-\tau\sqrt{1+c^2}\,\delta}\vert\mbox{\boldmath $b$}\vert^ne^R\left(\frac{1}{R^{n/2}}+1\right)
+
+\frac{\vert\mbox{\boldmath $b$}\vert^{2n}}
{R^n}e^{\displaystyle -(\vert\mbox{\boldmath $b$}\vert^2(\delta/c)^2(1+c^2))/(4R)}\right\}.
\tag {5.17}
$$

\noindent
Now Proposition 5.2 is a direct consequence of (5.12) and (5.17).

\noindent
$\Box$

{\bf\noindent Remark 5.2.}
All the derivatives of $K_z(x,t)$ also have a similar property:
for each $\alpha\in\mbox{\boldmath $Z$}_+^n$ and $\beta\in\mbox{\boldmath $Z$}_+$
$$\displaystyle
e^{\tau\sqrt{1+c^2}\delta}\sup_{\displaystyle (x\,\,t)^T\cdot\omega(c)<-\delta}\vert\partial^{\alpha}_x\partial^\beta_t K_z(x,t)\vert
$$
is at most algebraically growing as $\tau\longrightarrow\infty$.

As a corollary of Proposition 5.1 and remark 5.2 we obtain a Carleman type formula.

\proclaim{\noindent Corollary  5.1.}  Let $(x_0,t_0)\in\Omega\times]0,\,T[$ be an arbitrary fixed point.
Assume that $T>0$, $\omega$, $\Gamma$ and $U$ satisfy (4.8), (4.9) and (4.10).
Let $v(x,t)=K_z(x-x_0,t-t_0)$ for $z$ given by (4.6) and $u$ be a solution of (2.7).
Then we have
$$
\displaystyle
u(x_0,t_0)=-
\lim_{\tau\longrightarrow\infty}I(\tau),
\tag {5.18}
$$
where
$$\begin{array}{c}
\displaystyle
I(\tau)=
\int_{\Gamma}\left\{\left(\frac{\partial v}{\partial\nu}(x,t)+\rho(x)v(x,t)\right)u(x,t)
-h_0(x,t)v(x,t)\right\}dSdt
-\int_{U}v (x,0)u(x,0)dx.
\end{array}
$$

\endproclaim

{\bf\noindent Remark 5.3.} Yarmukhamedov \cite{Y} considered the
Cauchy problem for the Laplace equation in a three-dimensional
bounded domain $D$ that is bounded by the plane $x_3=0$ and by
smooth surfaces lying in the half-space $x_3>0$. He gave a formula
for calculating the value of the solution at a given point inside
the domain from the Cauchy data on the portion in $x_3>0$ of
$\partial D$. For the purpose he made use of a special fundamental
solution for the Laplace operator which has been introduced by
himself in \cite{Yo} and is parameterized by an entire function
$E(w)$, $w\in\mbox{\boldmath $C$}$ satisfying a suitable growth
condition when $\vert\text{Im}\,w\vert\longrightarrow\infty$. His
fundamental solution $\Phi_E$ takes the form for $x=(x_1,x_2,x_3)$
with $x'=(x_1,x_2)\not=0$:
$$\displaystyle
\Phi_E(x)
=-\frac{1}{2\pi^2}\int_0^{\infty}\text{Im}\,\left(\frac{E(x_3+i\sqrt{\vert x'\vert^2+u^2}\,)}{x_3+i\sqrt{\vert x'\vert^2+u^2}}\right)
\frac{du}
{\sqrt{\vert x'\vert^2+u^2}}.
$$
He chose the special $E(w)=e^{\tau w}$ with $\tau>0$.  Then $\Phi_E(x)$
has the representation
$$\displaystyle
\Phi_E(x)
=-\frac{e^{\tau x_3}}{2\pi^2}
\int_0^{\infty}
\left(x_3\frac{\sin\,(\tau\sqrt{\vert x'\vert^2+u^2}\,)}{\sqrt{\vert x'\vert^2+u^2}}
-\cos\,(\tau\sqrt{\vert x'\vert^2+u^2}\,)\right)\frac{du}{\vert x\vert^2+u^2}.
$$
From this one sees that the Cauchy data of $\Phi_E(y-x)$ on the portion in $y_3=0$ of
$\partial D$ for an arbitrary fixed $x\in D$ decays exponentially
as $\tau\longrightarrow\infty$.  This fact corresponds to
Proposition 5.2 and is an evidence that the formula (5.18) can be
considered as an extension to the heat equation of his formula.

It would be interesting to find a hidden `parameter' for $K_z$ like $E$ for $\Phi_E$.
That is: can one find a family of special fundamental solutions for the backward heat operator
$\partial_t+\triangle$ that contains $K_z$ as a special member?  In other words, can one find another
fundamental solution for the backward heat operator that is decaying in one side of a hyper surface
not plane?

$$\quad$$

\centerline{{\bf Acknowledgement}}

This research was partially supported by Grant-in-Aid for
Scientific Research (C)(No.  18540160) of Japan  Society for
the Promotion of Science.

$$\quad$$

\vskip1cm
\noindent
e-mail address

ikehata@math.sci.gunma-u.ac.jp

\end{document}